\documentclass[a4paper]{amsart}

\usepackage[margin=0.75in,headheight=15pt,footskip=25pt]{geometry}
 
\usepackage{amsmath}
\usepackage{amssymb}
\usepackage{amsthm}
\usepackage{subcaption}

\usepackage{tikz}
\usetikzlibrary{calc,patterns,angles,quotes}
 
\usepackage{enumerate}
\usepackage{cite} 
\usepackage{hyperref}

\hypersetup{colorlinks=true,allcolors=blue}
\usepackage[usenames,dvipsnames]{pstricks}
 
\usepackage{siunitx,booktabs}
\usepackage[title]{appendix}

\newtheorem{theorem}{Theorem}[section]

\newtheorem{lemma}[theorem]{Lemma}
\newtheorem{definition}{Definition}[section]
\newtheorem{remark}{Remark}[section]
\newtheorem{assumption}{Assumption}[section]

\numberwithin{figure}{section}
\numberwithin{equation}{section}
\newtheorem{algorithm}{Algorithm}

\newcommand{\tribar}{\vert\kern-0.25ex\vert\kern-0.25ex\vert}
\newcommand{\n}{\nu}

\newcommand{\M}{\bold{M}}
\renewcommand{\S}{\bold{S}}
\newcommand{\Q}{\bold{T}}
\newcommand{\D}{\bold{D}}
\newcommand{\C}{\mathcal{C}}

\renewcommand{\P}{\mathbb{P}}

\newcommand{\ddiv}{\text{div}\,}

\newcommand{\Gn}{\mathcal{G}^n}
\newcommand{\Aa}{\bold{A}_{\widehat\beb}^1}
\newcommand{\Ab}{\bold{A}^2}
\newcommand{\ba}{{\bold b}_{\widehat\alb,\widehat\beb}}
\newcommand{\bb}{{\bold b}^2}
\newcommand{\Aan}{\bold{A}_{\beb^n}^1}
\newcommand{\Abn}{\bold{A}^2}
\newcommand{\ban}{{\bold b}_{\alb^n,\beb^n}}

\newcommand{\Th}{\mathcal{T}_h}
\newcommand{\Sh}{\mathcal{S}_h}
\newcommand{\BM}{\mathcal{B}_M}
\newcommand{\Bn}{\mathcal{\widetilde {B}}_{n}}
\newcommand{\N}{N}
\newcommand{\Nh}{\mathcal{N}_h}
\newcommand{\Nhb}{\mathcal{N}_h^{\partial}}
\newcommand{\Nhi}{\mathcal{N}_h^{in}}
\newcommand{\Ehh}{\mathcal{E}_h}

\newcommand{\NT}{N_0}
\newcommand{\Zh}{\mathcal{Z}_h}
\newcommand{\Oh}{\mathcal{O}}

\newcommand{\Rn}{\mathbb{R}^{\N}}
\newcommand{\R}{\mathbb{R}}
\newcommand{\al}{{\alpha}}
\newcommand{\alb}{\boldsymbol{\alpha}}
\newcommand{\bchi}{\boldsymbol{\chi}}
\newcommand{\bpsi}{\boldsymbol{\psi}}
\newcommand{\be}{{\beta}}
\newcommand{\beb}{\boldsymbol{\beta}}

\newcommand{\bff}{\bold{f}}
\newcommand{\bw}{\bold{w}}
\newcommand{\q}{\bold{q}}
\newcommand{\bv}{\boldsymbol{v}}
\newcommand{\x}{\bold{x}}
\renewcommand{\n}{\bold{n}}

\title[Error analysis of positive preserving schemes for Chemotaxis]{Error analysis of a backward Euler positive preserving stabilized scheme for a Chemotaxis system}
\author{P. Chatzipantelidis}
\address{(1) Department of Mathematics and Applied Mathematics,
University of Crete, 71003, Heraklion, Greece}
\address{(2) Institute of Applied and Computational Mathematics, FORTH, Heraklion, 70013, Greece}
\email{p.chatzipa@uoc.gr}
\author{C. Pervolianakis}
\address{Institut für Mathematik, Friedrich-Schiller-Universität Jena, 07743, Jena, Germany}
\email{christos.pervolianakis@uni-jena.de}


\date{today is \today}


\begin{document}

\begin{abstract}
For a Keller-Segel model for chemotaxis in two spatial dimensions we consider a modification of a positivity preserving fully discrete scheme using a local extremum diminishing flux limiter. We discretize space using piecewise linear finite elements on an quasiuniform triangulation of acute type
and time by the backward Euler method. We assume that initial data are sufficiently small in order not to have a blow-up of the solution. Under appropriate assumptions on the regularity of the exact solution and the time step parameter we show existence of the fully discrete approximation and derive error bounds in $L^{2}$ for the cell density and $H^{1}$ for the chemical concentration. We also present numerical experiments to illustrate the theoretical results. 
\end{abstract}
\maketitle

\section{Introduction} \label{sec:Intro}
We shall consider 
a Keller-Segel system of equations of parabolic-parabolic type, where we seek
$u=u(\x,t)$ and $c=c(\x,t)$  for $(\x,t)\in{\Omega}\times [0,T],$ satisfying 
\begin{equation}\label{Minimal_model_uc}
\begin{cases}
u_t = \mu\Delta u - \lambda\,\ddiv( u\nabla c ), & \text{in }{{\Omega}}\times [0,T],\\
c_t = \eta\Delta c - c + u, &  \text{in }{{\Omega}}\times [0,T],\\
\frac{\partial u}{\partial \n}  =0,\;\;\;\frac{\partial c}{\partial \n} = 0,&\text{on }\partial {\Omega}\times[0,T],\\
u(\cdot,0)  = u^0,\quad c(\cdot,0)=c^0,&\text {in }{{\Omega}},
\end{cases}
\end{equation}
where ${\Omega}\subset{\R}^2$ is a convex bounded domain with boundary $\partial\Omega$, $\n$ is the outer unit 
normal vector to $\partial\Omega$, $\partial/\partial \n$ denotes differentiation along $\n$ on 
$\partial\Omega$, $\mu$, $\eta$, $\lambda $ are positive constants and $u^0, c^0\ge0$, $u^0\neq0$. 

The chemotaxis model \eqref{Minimal_model_uc}  describes the aggregation of slime molds resulting from their 
chemotactic features, cf. e.g. \cite{keller1970}. The function $u$  is the cell density of cellular slime molds, $c$ is the concentration of the chemical substance secreted by molds themselves, $c-u$ is the ratio of generation or extinction, and $\lambda$ is a chemotactic sensitivity constant. The parameters $\mu,\,\eta$ are positive constants that
represent, the cell diffusion and the chemical diffusion constant, respectively.

There exists an extensive mathematical study of chemotaxis models, cf. e.g., \cite{perthame2007,hillen2009,horstmann2003,horstmann2004,suzuki2005} and references therein. It is well-known that the solution of  \eqref{Minimal_model_uc} may blow up in finite time.  However, if $\|u_0\|_{L^{1}(\Omega)}\leq 4\pi\lambda ^{-1},$  the solution $(u,c)$ of \eqref{Minimal_model_uc} exists for all time and is bounded in $L^{\infty}$, cg. e.g. \cite{nagai1997}.

A key feature of the system \eqref{Minimal_model_uc} is the conservation of the solution $u$ in $L^{1}$ norm, $\|u(t)\|_{L^{1}(\Omega)} = \|u^0\|_{L^{1}(\Omega)}$, for $0\leq t\leq T$, which is an immediate result of the preservation of non-negativity of  $u$, and the conservation of total mass $\int_\Omega\,u(\x,t)\,d\x = \int_\Omega\,u_0(\x)\,d\x$, for $0\leq t\leq T$.

Capturing blowing up solutions numerically is a challenging problem and many numerical methods have been proposed
to address this. 
The main difficulty in constructing suitable numerical schemes is to preserve several essential properties of the
 Keller-Segel equations such as positivity, mass conservation, and energy dissipation.

Some numerical schemes were developed with positive-preserving conditions, cf. e.g. 
\cite{filbet2006,chertock2008,chertock2018}, which depend on a particular spatial discretization and impose 
CFL restrictions on the time step.
Other approaches include,  finite-volume based  numerical methods, \cite{chertock2008,chertock2018}, high-order 
discontinuous Galerkin methods, \cite{epshteyn2009,epshteyn2008,li-shu2017}, a flux corrected finite element method, 
\cite{strehl2010,strehl2013}, and a  novel numerical method based on symmetric reformulation of the chemotaxis system, \cite{liu2018}.
For a more detailed review on recent developments of numerical methods for chemotaxis problems, we refer to 
\cite{chertock2019,tadmor}.

In order to maintain the total mass and the non-negativity of the numerical approximations of the system 
\eqref{Minimal_model_uc}, Saito in \cite{saito2012, saito2007} proposed and analyzed a fully discrete method 
that uses an upwind finite element scheme in space and backward Euler method in time. His proposed finite element 
scheme made use of Baba and Tabata's upwind approximation, see \cite{tabata1981}. Strehl \textit{et al.} in \cite{strehl2010} 
proposed a slightly different approach. 
The stabilization was implemented  at a pure algebraic level via algebraic flux correction, see \cite{kuzmin2010}. 
This stabilization technique can be applied to piecewise linear finite element methods to counter the effect of the positive-preserving conditions in the convergence order and maintain the mass conservation and the non-negativity of the solution. 

For the forthcoming analysis and without loss of generality, we assume throughout the paper that $\mu, \eta = 1$. In the  variational form of \eqref{Minimal_model_uc}, we seek $u(\cdot,t)\in H^{1}$ and 
$c(\cdot,t)\in H^{1},$ for $t\in[0,T]$, such that
\begin{equation}\label{weak_uc}
\begin{aligned}
(u_t, v) + (\nabla u -\lambda u\nabla c, \nabla v) & =  0,\quad\forall v\in H^{1},\quad \text{ with }u(0)   = u^0,\\
(c_t, v) + (\nabla c, \nabla v) +(c-u, v) & = 0,\quad\forall\,v\in H^{1}, \quad \text{ with }c(0)   = c^0,
\end{aligned}\end{equation}
where $(f,g)=\int_\Omega fg\,dx$.

In our analysis we consider regular triangulations $\Th=\{K\}$ of $\overline\Omega$, with $h = \max_{K\in \Th}h_K$, $h_K=\text{diam}(K)$, 
and the finite element spaces 
\begin{equation*}
\Sh : = \{ \chi \in \C: \chi\vert_{K} \in \mathbb{P}_1,\forall K\in\Th\},
\end{equation*} 
where ${\C}={\C}(\overline{{\Omega}})$ denotes the continuous functions on $\overline{{\Omega}}$.

A semi-discrete approximation of the variational problem \eqref{weak_uc} 
is: Find $u_h(t)\in \Sh$ and $c_h(t)\in \Sh$, for $t\in[0,T]$,  with  $u_h(0)  = u_h^0\in\Sh$ and  $c_h(0)  = c_h^0\in\Sh$, such that
\begin{equation}\label{fem_uc}
\begin{aligned}
(u_{h,t},\chi) + (\nabla u_h - \lambda  u_h\nabla c_h, \nabla \chi) & = 0,\ \forall\chi\in \Sh,\\
(c_{h,t},\chi) + (\nabla c_h,\nabla \chi) + (c_h-u_h, \chi) & = 0,\ \forall\chi\in \Sh.
\end{aligned}
\end{equation}

We now formulate   \eqref{fem_uc} in matrix form.  Let  
 $\Zh = \lbrace Z_j\rbrace_{j=1}^{\N}$  be the set of nodes  in $\Th$ and $\lbrace \phi_j \rbrace_{j=1}^{\N}\subset \Sh$ the corresponding nodal basis, with $\phi_j(Z_i)=\delta_{ij}$.
Then, we may write
$u_h(t)=\sum_{j=1}^{\N}\al_j(t)\phi_j$, with $u_h^0=\sum_{j=1}^{\N}\al_j^0\phi_j$ and $c_h(t)= \sum_{j=1}^{\N}\be_j(t)\phi_j$, with $c_h^0=\sum_{j=1}^{\N}\be_j^0\phi_j$.
Thus, the semi-discrete problem \eqref{fem_uc} may then be expressed,   with 
$\alb = (\al_1, \dots, \al_{\N})^T$ and $\beb = (\be_1, \dots , \be_{\N})^T$, as
\begin{equation}\begin{aligned}
{\M}\alb^\prime(t) + (\S  - \Q_{\beb}) \alb(t) & = 0,\ \,\quad\qquad\text{ for } t\in[0,T], \text{ with }\alb(0)=\alb^0,\label{matrix_u}\\
{\M}\beb^\prime(t)+(\S +\M)\beb(t) & = \M\alb(t),\quad\text{ for } t\in[0,T], \text{ with } \beb(0)=\beb^0,
\end{aligned}\end{equation}
where $\alb^0 = (\al_1^0, \dots, \al_{\N}^0)^T$,  $\beb^0 = (\be_1^0, \dots , \be_{\N}^0)^T$, $\M=(m_{ij})$,  $m_{ij}=(\phi_i,\phi_j)$,  $\S=(s_{ij})$, 
 $s_{ij}=(\nabla\phi_i,\nabla\phi_j)$, 
\begin{align}\label{tau_def_matrix}
\Q_{\beb}=(\tau_{ij}(\beb)) \quad\text{and}\quad \tau_{ij}(\beb) =\lambda \sum_{\ell=1}^{\N}\be_\ell(\phi_j\nabla \phi_\ell,\nabla\phi_i), \quad \text{for } i,j=1,\dots,\N.
\end{align}
We will often suppress the index $\beb$ in the coefficients $\tau_{ij}=\tau_{ij}(\beb)$ and in $\Q=\Q_{\beb}$. An immediate property of the matrix $\Q_{\beb}$ is the zero column-sum for all $\beb\in\Rn$, that is
\begin{align}\label{zero_sum}
\sum_{i=1}^{\N} \tau_{ij}(\beb) = \lambda\sum_{\ell=1}^{\N}\be_\ell(\phi_j\nabla \phi_\ell, \nabla \sum_{i=1}^{\N}\phi_i) = 0,
\end{align}
since $\sum_{i=1}^{\N}\phi_i=1.$

 The matrices $\M$ and $\S$ are both symmetric 
 and positive definite, however   $\Q$ due to the chemotactical flux 
 $\lambda  u(t)\nabla c(t)$  is not symmetric.

Note that the semi-discrete solutions $u_h(t), c_h(t)$ of \eqref{fem_uc} are non-negative if and only
 if the coefficient vectors $\alb(t),\beb(t)$ are  non-negative element-wise. 
 In order to ensure non-negativity, we employ the lumped mass method, which results from replacing the 
 mass matrix $\M$ in \eqref{matrix_u} 
 with a diagonal matrix $\M_L$ with elements $\sum_{j=1}^{\N}m_{ij}$. 
 
 A sufficient condition for $\alb(t)$ to be  non-negative element-wise is that 
 the off diagonal elements of $\S - \Q$ are  non-positive. Further,  for $\beb(t)$ to be  non-negative element-wise, it suffices  that  the off diagonal elements of $\S$ are non-positive.  
 
Assuming, that $\Th$ satisfies an  acute condition, i.e., all interior angles of a triangle $K\in\Th$ are less or equal than $\pi/2$, we have that $s_{ij}\le0$, cf. e.g., \cite{draganescu2004}.  Then, in order to ensure that the 
off diagonal elements of $\S - \Q$ are  non-positive we may  add an artificial diffusion operator $\D=\D_{\beb}$. This technique is commonly used in conservation laws,  cf. e.g. \cite{kuzmin2010} and references therein. 
This modification of the semi-discrete scheme \eqref{weak_uc} is proposed in \cite{strehl2010}. This scheme is often called
 \textit{low-order scheme} since we introduce an error which manifests in the order of convergence.

To improve the convergence order of the \textit{low-order scheme}, Strehl \textit{et al.} in \cite{strehl2010} proposed another scheme, which is called \textit{algebraic flux correction scheme} or \textit{AFC scheme}. To derive the AFC scheme we decompose the error, introduced in the low-order scheme by adding the artificial diffusion operator, into internodal fluxes. Then we appropriately restore the optimal accuracy in regions where the solution does not violate the non-negativity. 

The AFC method gained a systematical attention in the last two decades and the first work that provide a rigorous mathematical proofs was established in \cite{gabriel2016}, in the context of the steady-state convection diffusion equation. In the latter work, the authors proved the solvability of the nonlinear scheme that arises from the interpretation of the AFC method as well as they derive error estimates in energy norm. 

In particular, they proved that, in convection--dominant regime, the AFC scheme convergences in energy norm with at least order $\mathcal{O}(h^{1/2}),$ which term is due to the stabilization term. 
By properly choosing the algorithm to limit the internodal fluxes so that is linearity preserving, it can improved the latter estimate to $\mathcal{O}(h),$ see, e.g., \cite{gabriel2018}.

Given this, we will consider limiters that satisfy the discrete maximum principle and preserve linearity on arbitrary meshes, such as the one proposed by Barrenechea \textit{et al.} \cite{gabriel2017b}. In our scheme, we will not limit the diffusion terms, since we work on triangulations consisting of weakly acute or right triangles, for which the stiffness matrix has the correct sign pattern. This allows us to derive an optimal estimate for the stabilization term introduced by the AFC method.

Our purpose here, is to analyze fully discrete schemes, for the approximation of \eqref{Minimal_model_uc}, 
by discretizing in time the \textit{low-order scheme} and the \textit{AFC scheme} using the backward  Euler method. 
We will consider the case where the solution of \eqref{Minimal_model_uc} remains bounded for all $t\ge0$, therefore  we will assume that $\|u_0\|_{L^{1}(\Omega)}\leq 4\lambda ^{-1}\pi.$ 

Our analysis of the stabilized schemes is based on the corresponding one employed by Barrenechea \textit{et al.} in \cite{gabriel2016}. In order to show existence of the solutions of the nonlinear fully discrete schemes, 
we employ a fixed point argument and demonstrate that our approximations remain uniformly bounded, provided that 
$\Th$ is quasiuniform and our time step parameter $k$, is such that $k=O(h^{1+\epsilon})$, with $0<\epsilon<1$.

We shall use standard notation for the Lebesgue  and Sobolev spaces, namely we denote  $W^{m}_p=W^{m}_p(\Omega)$, $H^{m}=W^{m}_{2}$, $L^{p}=L^{p}(\Omega)$, and  with $\|\cdot\|_{m}=\|\cdot\|_{H^{m}}$,  $\|\cdot\|_{L^{p}}=\|\cdot\|_{L^{p}(\Omega)}$, $\|\cdot\|=\|\cdot\|_{L^{2}}$for $m\in\mathbb{N}$ and $p\in[1,\infty]$, the corresponding norms. For the problem considered here, Saito in \cite{saito2012} for piecewise linear elements has derived suboptimal error estimates in $W^{1}_{\infty}$ and $L^{p}$ norm,  $p>2$, for a mass conservative finite element scheme and Epshteyn \textit{et al.} in \cite{epshteyn2009,epshteyn2008}, using a Discontinuous Galerkin method with piecewise quadratic elements,  also derived error estimates in $H^{1}$ norm of order $\Oh(h^2)$, assuming sufficiently smooth solutions.

The fully discrete schemes we consider approximate $(u^n, c^n)$ by $(U^n,C^n)\in \Sh\times\Sh$ where $u^n=u(\cdot,t^n)$, $c^n=c(\cdot,t^n)$, $t^n=nk$, $n=0,\dots, \NT$ and $\NT\in\mathbb{N}$, $\NT\ge1$, $k=T/\NT$. Assuming that the solutions $(u,c)$ of \eqref{Minimal_model_uc} are sufficiently smooth, with $u\in W^{1}_{\infty},\,c\in W^{2}_{\infty}$, we derive error estimates 
of the form 
 \begin{align*}
\|U^n-u^n\| & \le C ( k + k^{-1/2}h^2+h^{3/2}|\log h|),\quad
\|C^n-c^n\|_{1}  \le C ( k + k^{-1/2}h^2+h).
\end{align*}

The paper is organized as follows: In Section \ref{section:preliminaries} we introduce notation and the 
semi-discrete  \textit{low-order scheme} and the \textit{AFC scheme} for the discretization of \eqref{Minimal_model_uc}. Further, we prove some auxiliary results for  the stabilization terms, that we will employ in the analysis that follows and 
rewrite  the low-order and AFC scheme, as general semi-discrete scheme. 
In Section \ref{section:fully_discrete}, we discretize the  general semi-discrete scheme
in time, using the backward  Euler method. For a sufficiently smooth solution of  \eqref{Minimal_model_uc} and $k=O(h^{1+\epsilon})$, with $0<\epsilon<1$, we demonstrate that there exists a unique discrete solution which remains bounded and 
derive error estimates in $L^{2}$ for the cell density and $H^{1}$ for the chemical concentration. In Section \ref{section:positivity}, we show that the discrete solution preserves the positivity. Finally, in Section \ref{section:numerical_results}, we present numerical experiments, illustrating our theoretical results.

\section{Preliminaries}\label{section:preliminaries}

We consider a  family of regular triangulations $\Th=\{K\}$ of a convex polygonal domain $\overline\Omega\subset\R^2$.
We will assume that the family $\Th$ satisfies the following assumption.

\begin{assumption}\label{mesh-assumption}
Let $\Th=\{K\}$ be a  family of regular triangulations  of $\overline\Omega$ such that any edge of any $K$ is either a subset of the boundary $\partial\Omega$ or an edge of another $K \in \Th$, and in addition
\begin{enumerate}
\item  $\Th$ is shape regular, i.e, 
there exists a constant $\gamma>0,$ independent of $K$ and $\Th,$ such that 
\begin{equation}\label{shape_regularity}
\frac{h_K}{\varrho_K} \leq \gamma,\quad \forall K\in\Th,
\end{equation}
where $\varrho_K=\text{diam}(B_K)$, and $B_K$ is the inscribed ball in $K$.
\item The  family of triangulations $\Th$ is quasiuniform, i.e., there exists constant $\varrho>0$ such that
\begin{align}\label{quasi-uniformity}
\frac{\max_{K\in\Th}h_K}{\min_{K\in\Th}h_K} \leq \varrho,\quad\forall K\in\Th,
\end{align}
\item All interior angles of $K\in\Th$ are less or equal than $\pi/2$.
\end{enumerate}
\end{assumption}

\begin{figure}
\centering
\includegraphics[scale=0.7]{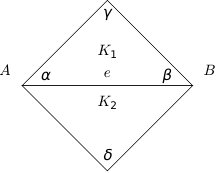} \hspace{.1cm}
\includegraphics[scale=0.7]{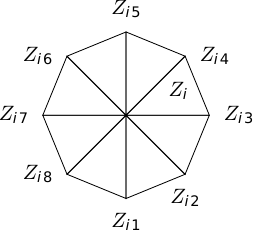} \hspace{.1cm}
{%
\begin{tikzpicture}[scale=0.3]
\tikzstyle{every node}=[font=\Large]
\draw  (3.75,5.25) -- (6.25,7.75);
\draw  (6.25,7.75) -- (5,10.25);
\draw [dashed] (6.25,7.75) -- (7.5,5.25);
\draw [dashed] (6.25,7.75) -- (8.75,10.25);
\draw [dashed] (8.75,10.25) -- (7.5,5.25);
\draw  (5,10.25) -- (3.75,5.25);
\node  [font=\scriptsize] at (6.7,9) {$Z_i$};
\node [font=\scriptsize] at (5.25,8) {$K_i$};
\node [font=\scriptsize] at (5,8) {$$};
\node [font=\scriptsize] at (7.5,8) {$\tilde K_i$};
\node [font=\scriptsize] at (4.5,10.5) {$Z_j$};
\node [font=\scriptsize] at (3,5) {$Z_{\ell}$};
\draw  (6.25,7.75) -- (5,5.25);
\draw   (6.25,7.75) -- (6.25,10.25);
\draw (6.25,10.25) -- (5,11.75);
\draw   (5,5.25) -- (3.5,4.);
\draw  (5,10.25) -- (6.25,10.25);
\draw   (3.75,5.25) -- (5,5.25);
\node [font=\scriptsize] at (3.5,8) {$\omega_i$};
\node [font=\scriptsize] at (3,6.5) {$\Omega$};
\end{tikzpicture}
}%
\caption{Left and Center:  Sub-domains, $\omega_e$ and $\omega_I$ of the triangulation $\Th.$ Right: A triangle $K_i$ and a symmetric extension $\tilde K_i$ at a boundary node $Z_i$.}\label{fig:patches}
\end{figure}

Let  $\Nh:=\{ i: Z_i \text{ a node of the triangulation } \Th\}$, $\Ehh$ be the set of all edges of the triangulation $ \Th$. 
We denote  $\omega_e$ the collection of triangles with a common edge $e\in\Ehh$,  see Fig. \ref{fig:patches},
and $\omega_i$, $i\in \Nh$, the collection of triangles with a common vertex $Z_i$,  
i.e. $\omega_i = \cup_{Z_i\in K}\overline{K},$ see Fig. \ref{fig:patches}.
The sets $\Zh(\omega)$ and $\Ehh(\omega)$ contain the vertices and  the edges, respectively, of a 
 subset of $\omega\subset\Th$, $\Nh(\omega)=\{ j: Z_j\in \Zh(\omega)\}$ and  
 $\Nh^i:=\{ j: Z_j\in \Zh, \text{adjacent to }Z_i\}$. Using the fact that $\Th$ is 
shape regular, there exists a constant $\kappa_\gamma$, independent of $h$,  such that the number of vertices in $\Nh^i$  is 
less than  $\kappa_\gamma$, for $i=1,\dots, \N$. 
Let  $\Nhb:=\{ i: Z_i\in\Zh \cap\partial\Omega\}$ and $\Nhi:=\Nh\setminus\Nhb$. In the sequel we will also consider additional "ghost" nodes outside $\Omega$. 
For $i\in \Nhb$ and $Z_i\in K\cap\partial\Omega$, $K_i\in\Th$, let $\tilde K_i$,  be the  symmetric extension of $K_i$, with respect to $Z_i$, which lies outside $\Omega$, $\tilde K_i\not\subset \Omega$, and $\tilde \omega_i=\cup_{Z_i\in K_i}{K_i\cup \tilde K_i}$, see Fig \ref{fig:patches}. 
Next, for $v\in \P_1(K_i)$, $i\in \Nhb$, let $L_{K_i}v\in \P_1(K_i\cup\tilde K_i)$ be the linear 
extension of $v$ in $K_i\cup\tilde K_i$   and for $v$ a piecewise linear function on $\omega_i$, $i\in\Nhb$, let $L_iv$ be the piecewise linear function on $\tilde\omega_i$, with $L_iv|_{K_i\cup\tilde K_i}=L_{K_i}v$ for $K_i\in\omega_i$.

Since $\Th$ satisfies \eqref{quasi-uniformity}, we have for all $\chi\in \Sh,$ cf., e.g., \cite[Chapter 4]{brenner2008}, 
\begin{equation}\label{eq:inverse_estimate}
\|\chi\|_{L^{\infty}} + \|\nabla \chi\|  \le Ch^{-1}\|\chi\|
\text{ and } \|\nabla \chi\|_{L^{\infty}}  \le Ch^{-1}\|\nabla\chi\|.
\end{equation} 
Further, in our analysis, we will employ the following trace inequality which holds for $K\in\Th,$ see, e.g., \cite[Theorem 1.5.1.10]{grisvard}, that for every triangle $K\in\Th,$ the following trace inequality holds,
\begin{equation}\label{eq:scaled_trace}
\|v\|_{L^{p}(\partial K)}  \leq C\left(h^{-1/p}_K\| v\|_{L^{p}(K)} + h_K^{1-1/p}\|\nabla v\|_{L^{p}(K)}\right),\;\;\forall\,v\in W^{1}_p(K),\;\;1\leq p\leq \infty.
\end{equation} 

We consider the $L^2$ projection $P_h:L^2\to\Sh$ and the elliptic projection $R_h:H^{1}\to \Sh$ defined by
\begin{align}
( P_hv - v, \chi) &= 0, \quad\forall \chi\in \Sh,\label{L2_projection_2D}\\
(\nabla R_hv - \nabla v, \nabla \chi) + (R_hv-v, \chi) &= 0,\quad\forall \chi\in \Sh\label{ritz_projection_2D}.
\end{align} 
In view of the mesh Assumption \ref{mesh-assumption}, $P_h$ and $R_h$ satisfy the following bounds, cf. e.g., \cite[Chapter 8]{brenner2008} and \cite{saito2012}. 
\begin{alignat}{2}
\|P_hv\|_{1} &\le C\|v\|_{1}, &\quad &\forall v\in H^{1},\label{L2_projection_stab}\\
\|v - R_hv\| + h\|v - R_hv\|_{1} &\le Ch^2\|v\|_{2}, &\quad &\forall v\in H^{2},\label{ritz_projection_est2_2D}\\
\|v - R_hv\|_{L^{\infty}} & \le Ch^2|\log h|\|v\|_{W^{2}_\infty}, &\quad &\forall v\in W^{2}_{\infty},\label{ritz_projection_inf}\\
\|\nabla R_h v\|_{L^{\infty}} &\le C\|v\|_{W^{1}_\infty}, &\quad &\forall v\in W^{1}_\infty.\label{ritz_projection_stab}
\end{alignat}

In the sequel we will present two stabilized semi-discrete schemes for the numerical approximation of \eqref{Minimal_model_uc}, namely the low-order scheme and the AFC scheme, which have been proposed in \cite{strehl2010}. 

\subsection{Low-order scheme}
The semi-discrete problem \eqref{fem_uc} may be expressed in the matrix form \eqref{matrix_u}, 
where the matrices $\M$ and  $\S$ are symmetric  and positive definite.
However $\Q_{\beb}$ is not symmetric, but as we showed, see, \eqref{zero_sum}, it has  a zero-column sum for all $\beb\in\Rn.$ 

For a function $v\in \C$, let $\bv=(v_1,\dots,v_\N)^T$ denote the vector with coefficients its nodal values $v_i = v(Z_i)$, $Z_i\in\Zh$, $i=1,\ldots,\N$. We will often express the coefficients $\tau_{ij}$, $i,j=1,\dots,\N$, of $\Q_{\beb}$ as  functions of an element 
$\psi\in\Sh$, $\tau_{ij}=\tau_{ij}(\bpsi)=\tau_{ij}(\psi)$,  such that 
$\psi=\sum_j\psi_j\phi_j\in \Sh$ and  $\bpsi= (\psi_1, \dots , \psi_{\N})^T$. Thus the elements  of $\Q_{\beb}=\Q=(\tau_{ij})$, 
  may be expressed equivallently as, 
 \begin{equation}\label{T_def}
 \tau_{ij}=\tau_{ij}(\bpsi)=\tau_{ij}(\psi) = \lambda (\phi_j\nabla \psi,\nabla\phi_i)=\lambda \sum_{\ell=1}^{\N}\psi_\ell(\phi_j\nabla \phi_\ell,\nabla\phi_i).
 \end{equation}

In order to preserve non-negativity of $\alb(t)$ and $\beb(t)$, a low order semi-discrete scheme of minimal model has been proposed, cf. e.g. \cite{strehl2010}, where  $\M$ is replaced by the corresponding lumped mass matrix $\M_L$ and an artificial artificial diffusion operator $\D=\D_{\beb}=(d_{ij})$ is added to $\Q$, to 
eliminate all negative off-diagonal elements of $\Q$, so that $\Q+\D\ge0$, element-wise.
Thus, assuming that $\Th$ satisfies an  acute condition, i.e., all interior angles of a triangle 
$K\in\Th$ are less or equal than $\pi/2$, gives  that $s_{ij}\le0$ and hence, the off diagonal elements of 
$\S - \Q-\D$ are non-positive, $s_{ij}-\tau_{ij}-d_{ij}\le0$, $i\neq j$, $i,j=1,\dots, \N$. 
However, note that assuming $\Th$ to be acute is not a necessary condition to preserve non-negativity of $\alb(t)$ or $\beb(t)$.
 
Also, we will often suppress the index $\beb$ in the coefficients $d_{ij}=d_{ij}(\beb)$, $i,j=1,\dots,\N$, or express them as  functions of an element $\psi\in\Sh$, $d_{ij}(\psi)=d_{ij}(\bpsi)$, such that $\psi=\sum_j\psi_j\phi_j\in \Sh$ and  $\bpsi = (\psi_1, \dots , \psi_{\N})^T$. Since, we would like our scheme to maintain  the mass,  $\D$ must be symmetric with zero row and column sums, cf. \cite{strehl2010},  which is true if  $\D=(d_{ij})_{i,j=1}^{\N}$ is defined by
\begin{equation}\label{D_def}
d_{ij} : = \max\{ -\tau_{ij},0, -\tau_{ji}\}=d_{ji}\ge0,\quad\forall j\neq i\ 
\text{ and }\ 
d_{ii} : = -\sum_{j\neq i}d_{ij}.
\end{equation}
Thus the resulting system for the approximation of \eqref{Minimal_model_uc} is expressed as follows, we seek $\alb(t),\beb(t)\in\Rn$ such that, for $t\in[0,T]$,
\begin{equation}\begin{aligned}
\M_L\alb^\prime(t) + (\S  - \Q_{\beb}-\D_{\beb}) \alb(t) & = 0,\quad\qquad\text{ with }\alb(0)=\alb^0,
\label{low_matrix_u}\\
\M_L\beb^\prime(t)+(\S  + \M_L)\beb(t) & = \M_L\alb(t), \text{ with } \beb(0)=\beb^0.
\end{aligned}\end{equation}

 Let for $w\in{\Sh}$, $d_h(w;\cdot,\cdot):{\C}\times {\C}\to{\R},$  be a bilinear form defined by
\begin{equation}\label{stab_term}
d_h(w;v,z) := \sum_{i,j=1}^{\N}\,d_{ij}(w)(v_i - v_j)z_i,\quad\forall v,z\in{\C},
\end{equation}
and  $(\cdot,\cdot)_h$ be an inner product in $\Sh$ that approximates $(\cdot,\cdot)$ and is defined
by
\begin{equation}\label{quadrature}
(\psi,\chi)_h = \sum_{K\in\Th}Q_h^K(\psi\chi),\ \text{ with }
Q_h^K(g) = \frac{1}{3}|K|\sum_{Z\in \Zh(K)}g(Z)\approx \int_K g(x)\,dx,
\end{equation}
with $\Zh(K)$ the vertices of $K\in\Th$ and $|K|$ the area of $K\in\Th$.
Then following \cite{gabriel2016},  the coupled system \eqref{low_matrix_u} 
can be rewritten as: 
Find $u_h(t), c_h(t)\in \Sh$, with $u_h(0)   = u_h^0\in\Sh$ 
and $c_h(0)   = c_h^0\in\Sh$, such that
\begin{equation}
\begin{aligned}
(u_{h,t},\chi)_h + (\nabla u_h - \lambda  u_h\nabla c_h, \nabla \chi) + d_h(c_h;u_h,\chi)  & = 0,\quad \forall\,\chi\in \Sh,
\label{low_fem_u_2D}\\
(c_{h,t},\chi)_h + ( \nabla c_h,\nabla \chi) + (c_h - u_h, \chi)_h & = 0,\quad \forall\,\chi\in \Sh.
\end{aligned}
\end{equation}
 
  We can easily see that $(\cdot,\cdot)_h$ induces  an equivalent norm to $\|\cdot\|$ on $\Sh$.
Thus,  there exist constants $C,\,C'$ independed on $h$, such that
\begin{equation}\label{mass_lump_equivalence}
C\|\chi\|_h  \leq \|\chi\|  \leq C'\|\chi\|_{h},\ \text{ with }
 \|\chi\|_h = (\chi, \chi)_h^{1/2},\quad\forall\chi\in \Sh.
\end{equation}

\subsection{Algebraic flux correction scheme}
 The replacement of the standard FEM discretization \eqref{matrix_u} 
 by the low-order scheme \eqref{low_matrix_u} 
 ensures non-negativity but introduces an error 
 which manifests the order of convergence, cf. e.g. \cite{strehl2013,kuzmin2010}. 
 Thus, following \textit{Strehl et al.} \cite{strehl2013}, one may ``correct" the semi-discrete scheme 
 \eqref{low_matrix_u} 
 by introducing a  flux correction term. Hence, 
 we also consider an algebraic flux correction (AFC) scheme, which involves the decomposition of this error into 
 internodal fluxes, which can be used to restore high accuracy in regions where the solution is well resolved and no 
 modifications of the standard FEM are required. There exists various algorithms to implement
an AFC scheme. Here we will follow the one proposed by Barrenachea \textit{et. al.} in \cite{gabriel2017b}.

The AFC scheme is constructed in the following way.
Let $\bff=(f_1,\dots,f_\N)^T$  denote  the error of inserting the operator $\D_{\beb}$ in \eqref{matrix_u}, i.e., $\bff(\alb,\beb) =  - \D_{\beb}\alb.$ Using the zero row sum property of matrix $\D_{\beb}$, cf. \eqref{D_def}, we can show that the residual admits a conservative decomposition into internodal fluxes,
\begin{equation*}
f_i = \sum_{j\neq i}f_{ij},\quad f_{ji} = -f_{ij},\quad i = 1,\dots,\N,
\end{equation*}
where the amount of mass transported by the raw \textit{antidiffusive flux} $f_{ij}$ is given by
\begin{equation*}
f_{ij} := f_{ij}(\alb, \beb) = (\al_i - \al_j)d_{ij}(\beb),\quad\forall j\neq i,\quad i,j=1,\dots,\N.
\end{equation*}
For the rest of this paper we will call the internodal fluxes as anti-diffusive fluxes. Some of these anti-diffusive fluxes are harmless but others may be responsible for the violation of non-negativity. Such fluxes need to be canceled or limited so as to keep the scheme non-negative. Thus, every anti-diffusive flux $f_{ij}$ is multiplied by a solution-depended correction factor $\mathfrak{a}_{ij}\in[0,1]$, to be defined in the sequel, before it is inserted into the equation. Hence, the AFC scheme is the following: We seek $\alb(t),\beb(t)\in\Rn$ such that, for $t\in[0,T]$,
\begin{equation}\begin{aligned}
\M_L\alb^\prime(t) +(\S  - \Q_{\beb}-\D_{\beb})\alb(t)& = \overline{\bff}(\alb(t),\beb(t)),\quad \text{ with }\alb(0)=\alb^0,\label{ode_u_afc_2D}\\
\M_L\beb^\prime(t) +(\S  + \M_L)\beb(t)& =  \M_L\alb(t),\quad\qquad\text{ with }\beb(0)=\beb^0,
\end{aligned}\end{equation}
where $\overline{\bff}(\alb(t),\beb(t))=(\overline{f}_{1},\dots, \overline{f}_{\N})^T$, with
\begin{equation*}
\overline{f}_{i}:=\overline{f}_{i}(\alb(t),\beb(t))=\sum_{j\neq i}\mathfrak{a}_{ij}f_{ij},
\quad i=1,\dots,\N,
\end{equation*}
and $\mathfrak{a}_{ij}\in[0,1]$ are appropriately defined in view of the antidiffusive fluxes $f_{ij}$.

In order to determine the coefficients $\mathfrak{a}_{ij}$, one has to fix  a set of 
non-negative  coefficients $q_{i}$. In principle the choice  
of these parameters $q_{i}$ can be arbitrary. But efficiency and accuracy can dictate a strategy, which does not depend 
on the fluxes $f_{ij}$ but on the type of problem ones tries to solve and the mesh parameters. We will not elaborate more 
on the choice of $q_{i}$, for a more detail presentation we refer to \cite{gabriel2016,kuzmin2010}  and the references 
therein. In the sequel we
will employ two particular choices of $q_i$, cf. Lemma \ref{corrollary:linearity_preservation_limiters}. 

Due to the Neumann boundary conditions, we will need to modify the algorithm for determining  the coefficients $\mathfrak{a}_{ij}$ in \cite{gabriel2016,kuzmin2010}, cf.  \cite[Remark 10.45]{barrenechea2025}. 

For $v\in \P_1(K_i)$, $i\in \Nhb$, let $L_{K_i}v\in \P_1(K_i\cup\tilde K_i)$ be the linear 
extension of $v$ in $K_i\cup\tilde K_i$   and for $v$ a piecewise linear function on $\omega_i$, $i\in\Nhb$, let $L_iv$ be the 
piecewise linear function on $\tilde\omega_i$, with $L_iv|_{K_i\cup\tilde K_i}=L_{K_i}v$ for $K_i\in\omega_i$.

\begin{definition}\label{def-max} Let $\alb=(\al_1,\dots,\al_\N)^T$. Then for $i\in\Nhi$, $\al_i^{\max}$, $\al_i^{\min}$ are the local maximum and local minimum of  $\al_j$, 
$j\in\Nh(\omega_i)$. For $i\in \Nhb$, let  $v_{\al_i}$
 be the piecewise linear function on $\omega_i$, with $v_{\al_i}(Z_j)=\al_j$, $j\in\Nh(\omega_i)$.   Then $\al_i^{\max}$, $\al_i^{\min}$, $i\in\Nhb$, are the local maximum 
 and local minimum of  $L_iv_{\al_i}$ on $\tilde\omega_i$.
\end{definition}

To ensure that the AFC scheme maintains the non-negativity property, it is sufficient to choose the correction factors 
$\mathfrak{a}_{ij}$ such that the sum of anti-diffusive fluxes is constrained, cf. e.g., \cite{kuzmin2010}, as follows.
Let $q_i>0$, $i\in\Nh$, are given constants that do not depend on $\alb$ and
\begin{equation*}
{Q}_{i}^{+} = q_i(\al_i^{\max} - \al_i)\quad
\text{ and }\quad{Q}_{i}^{-} = q_i(\al_i^{\min} - \al_i),\quad i\in\Nh,
\end{equation*}
with $\al_i^{\max}$, $\al_i^{\min}$ given by Definition \ref{def-max}. The constants $\mathfrak{a}_{ij}$
should satisfy 

\begin{equation}\label{led_2D}
{Q}_{i}^{-} 
\leq \sum_{j\neq i}\mathfrak{a}_{ij}f_{ij}
\leq {Q}_{i}^{+},\quad i\in\Nh.
\end{equation}
\begin{remark}
Note that if all the correction factors $\mathfrak{a}_{ij}=0$, then the AFC scheme 
\eqref{ode_u_afc_2D}   
reduces to the low-order scheme \eqref{low_matrix_u}.  
\end{remark}
\begin{remark}\label{remark:led}
The criterion \eqref{led_2D} by which the correction factors are chosen, implies that the limiters used in \eqref{ode_u_afc_2D}  
guarantee that the scheme is non-negative. In fact, if $\al_i=\al_i^{\max}$ then \eqref{led_2D} implies the 
cancellation of all positive fluxes. Similarly, all negative fluxes are canceled if $\al_i=\al_i^{\min}$. In other words, 
a local maximum cannot increase and a local minimum cannot decrease. 
As a consequence, $\mathfrak{a}_{ij}f_{ij}$ cannot create an undershoot or overshoot at node $i.$ 
\end{remark}

We shall compute the correction factors $\mathfrak{a}_{ij}$ using Algorithm  \ref{algorithm-1},
 which is based on  Kuzmin, cf. \cite[Section 4]{kuzmin2010}. Then we have that
\begin{equation*}
Q_{i}^- \leq  R_{i}^-P_{i}^- \leq \sum_{j\neq i}\mathfrak{a}_{ij}f_{ij}\leq R_{i}^+P_{i}^+ \leq Q_{i}^+,\quad i\in\Nh,
\end{equation*}
which implies  that \eqref{led_2D} holds.

\noindent
\begin{algorithm}[Computation of correction factors $\mathfrak{a}_{ij}$] \label{algorithm-1}
Given data: 
\begin{enumerate}
\item The positive  coefficients $q_{i}$, such that $q_i = \Oh(h)$, $i\in\Nh$. 
\item The fluxes $f_{ij}$, $i\neq j$, $i,j=1,\dots,\N$.
\item The coefficients $\al_j,\,\be_j$,  $j=1,\dots,\N$.
\end{enumerate}

\noindent
Computation of factors $\mathfrak{a}_{ij}$.
\begin{enumerate}
\item
Compute the limited sums $P_{i}^{\pm}:= P_{i}^{\pm}(\alb,\beb)$, $i\in\Nh$, of positive and negative anti-diffusive fluxes
\begin{align*}
P_{i}^{+} = \sum_{j\neq i}\max\{0, f_{ij}\}
\quad \text{ and }\quad P_{i}^{-} = \sum_{j\neq i}\min\{0, f_{ij}\}, \quad i\in\Nh.
\end{align*}
\item
Retrieve the local extremum diminishing upper and lower bounds  ${Q}_{i}^{\pm} : = {Q}_{i}^{\pm}(\alb),\;i\in\Nh$,
\begin{align*}
{Q}_{i}^{+} = q_i(\al_i^{\max} - \al_i),
\quad\text{ and }\quad{Q}_{i}^{-} = q_i(\al_i^{\min} - \al_i),\quad i\in\Nh,
\end{align*}
where $\al_i^{\max}$, $\al_i^{\min}$ are the local maximum and local minimum given by Definition \ref{def-max}.
\item
Compute the coefficients $\overline{\mathfrak{a}}_{ij},$ for $j\neq i$, $i,j\in\Nh$,  which are defined by
\begin{equation}\label{correction_factors_definition}
\begin{aligned}
R_{i}^+=\min\left\lbrace 1, \frac{Q_{i}^+}{P_{i}^+}\right\rbrace,\quad R_{i}^-=\min\left\lbrace 1,\frac{Q_{i}^-}{P_{i}^-}\right\rbrace  \quad \text{and}\quad \overline{\mathfrak{a}}_{ij} = 
\begin{cases}
 R_{i}^+, &\text{if} \quad f_{ij} > 0,\\
1, &\text{if} \quad f_{ij} = 0,\\
 R_{i}^-, &\text{if} \quad f_{ij} < 0.
\end{cases}
\end{aligned}
\end{equation}
Note that if $P_{i}^-=0$ or $P_{i}^+=0$, then we define $R_{i}^-=1$ or $R_{i}^+=1$, respectively.
\item Finally,  the requested coefficients $\mathfrak{a}_{ij},$ for $j\neq i$, $i,j\in\Nh$  are defined by
\begin{align*}
\mathfrak{a}_{ij} :=\min\{\overline{\mathfrak{a}}_{ij},\overline{\mathfrak{a}}_{ji}\}, \text{ for } i,j\in\Nh.
%
\end{align*}
Note that then $\mathfrak{a}_{ij} = \mathfrak{a}_{ji}$, for $i,j\in\Nh$. 
\end{enumerate}  

\end{algorithm}

\begin{remark}
Following the definition of $\tau_{ij}$ in \eqref{T_def}, we may express $\mathfrak{a}_{ij}=\mathfrak{a}_{ij}(\bv,\bw)=\mathfrak{a}_{ij}(v,w)$, with $v=\sum_jv_j\phi_j$, $w=\sum_jw_j\phi_j$ and $\bv=(v_1,\dots,v_{\N})^T, \bw=(w_1,\dots,w_{\N})^T\in\Rn$.
\end{remark}

\begin{remark}\label{remark:linearity_preservation}
There exist  $\gamma_i\in\R$, $i\in\Nh$,  cf. \cite[Section 6]{gabriel2017b}, such that
\begin{equation}\label{def:gamma_i}
v_i - v_i^{\min} \leq \gamma_i (v_i^{\max} - v_i), \quad \forall v\in \P_1(\R^2),
\end{equation}
for $v_i=v(Z_i)$ and $v_i^{\max}$ and $v_i^{\min}$ the local maximum and local minimum given by Definition \ref{def-max} 
and
$$
\gamma_i=\dfrac{\max_{Z_j\in \partial\omega_i}|Z_i-Z_j|}{\mathrm{dist}(Z_i,\partial \omega_i^{\mathrm{conv}})},
$$
with $\omega_i^{\mathrm{conv}}$ the convex hull of $\omega_i$, for $i\in \Nhi$ and $\tilde\omega_i$ for $i\in \Nhb$. We note that, since we have assumed a weakly acute condition for $\Th,$ we have $\omega_i = \omega_i^{\mathrm{conv}}, \,i\in\Nhi.$ Further, if $\omega_i$ is symmetric with respect to $Z_i$ then $\gamma_i = 1$, see \cite{gabriel2017b}.
\end{remark}

\begin{definition}

The limiter $\overline{\mathfrak{a}}_{ij}$ defined in \eqref{correction_factors_definition} has the
 \emph{linearity preservation property} if
\begin{equation}\label{eqn:linear_preserve}
\overline{\mathfrak{a}}_{ij}(v,w) = 1,\quad \text{for }i,j\in\Nh,\text{ and }v\in {\P}_1(\R^2),\ w\in \Sh.
\end{equation} 
\end{definition}

\begin{lemma}\label{corrollary:linearity_preservation_limiters}
Let  the positive  coefficients $q_i,$ $i\in\Nh$, in Algorithm \ref{algorithm-1} be defined by
\begin{align}\label{def:q_i-1}
q_i : = \gamma_i \sum_{j\neq i} d_{ij},\quad i\in\Nh,
\end{align}
with $\gamma_i$ defined in  \eqref{def:gamma_i}, then the linearity preservation property
 \eqref{eqn:linear_preserve} is satisfied. Further, if there exists $M>0$ such that $\|\nabla w\|_{L^{\infty}} \leq M$ for $w\in\Sh$ and the constants   $q_i,$ $i\in\Nh$, in Algorithm \ref{algorithm-1},
are defined by
\begin{align}\label{def:q_i-2}
q_i : =  \gamma_i\frac{m_i}{\nu} \text{ with } \nu\in (0,1)\text{ and }\nu = \Oh(h^{1+\epsilon}),\ \epsilon\in (0,1),
\end{align}
with  $m_i$ the diagonal elements of $\M_L$,
 then there exists $h_M>0$ such that for $h\le h_M$, 
 \eqref{eqn:linear_preserve} is satisfied.
\end{lemma}
\begin{proof} 
The proof is based on the fact that the linearity preservation is equivalent to have
\begin{align}\label{LP_criterion}
Q_{i}^+ > P_{i}^+\ \text{ if } f_{ij} > 0\ \text{ and } Q_{i}^- < P_{i}^-\ \text{ if } f_{ij} < 0,
\end{align}
see \cite[Section 6]{gabriel2017b}. Let $v\in {\P}_1(\R^2)$ and
  the positive  coefficients $q_i,$ $i\in\Nh$, in Algorithm \ref{algorithm-1} are defined by \eqref{def:q_i-1}.
Following   the proof of \cite[Theorem 6.1]{gabriel2017b}  we can show \eqref{LP_criterion} and thus 
\eqref{eqn:linear_preserve}. Indeed using the definition of  $\gamma_i$ in Remark \ref{remark:linearity_preservation} and following \cite[Lemma 7]{gabriel2018} we obtain
\begin{align}
P_i^+ &= \sum_{j\neq i}\max\{0, f_{ij}\}  = \sum_{j\in \Nh^i}\max\{0, d_{ij}(w)(v_i - v_j)\} 
 = \sum_{j\in \Nh^i,\;v_i > v_j}d_{ij}(w)(v_i - v_j) \notag\\
& \leq \sum_{j\in \Nh^i,\;v_i > v_j}d_{ij}(w)(v_i - v_i^{\min})
 \leq \sum_{j\in \Nh^i}d_{ij}(w)(v_i - v_i^{\min})\notag\\
& \leq \gamma_i(v_i^{\max} - v_i)\sum_{j\in \Nh^i}d_{ij}(w) =q_i(v_i^{\max} - v_i)=Q_i^+.\label{ineq:P+}
\end{align}
Similarly, we can show $Q_{i}^- < P_i^-$, for $i\in\Nh$, cf. \cite[Section 6]{gabriel2017b} and hence \eqref{eqn:linear_preserve} is satisfied.
Let now  the positive  coefficients $q_i,$ $i\in\Nh$, in Algorithm \ref{algorithm-1} be defined by \eqref{def:q_i-2}.
For $w\in \Sh$, in view of the definition of $d_{ij}=d_{ij}(w)$ in \eqref{D_def}, 
\eqref{T_def} and the fact that the triangulation $\Th$ is shape regular, i.e., \eqref{shape_regularity},  there exists a constant $ C_\gamma>0$ independent of $h$, such that
\begin{align}
\vert d_{ij}(w)\vert & \leq \vert\tau_{ij}(w)\vert + \vert\tau_{ji}(w)\vert 
 \leq \|\nabla w\|_{L^{\infty}}\sum_{K\in \omega_i} (\|\nabla \phi_i\|_{L^{2}(K)}\|\phi_j\|_{L^{2}(K)} + \|\nabla \phi_j\|_{L^{2}(K)}\|\phi_i\|_{L^{2}(K)})\nonumber\\
& \leq C \|\nabla w\|_{L^{\infty}} \sum_{K\in \omega_i}h_K \leq C_\gamma h\|\nabla w\|_{L^{\infty}}.\label{est_dij_2D_L_infty}
\end{align}

Using now \eqref{est_dij_2D_L_infty} and the fact that the family of triangulations are shape regular, $m_i=\Oh(h^2)$, 
there exists $h_M>0$ such that for sufficiently small $h<h_M$, we get

\begin{align}
\sum_{j\in \Nh^i}d_{ij}(w) \leq C\,h\|\nabla w\|_{L^{\infty}} & = CM\,h\frac{m_i}{\nu}h^{\epsilon}\frac{\nu}{m_i}h^{-\epsilon} \leq CM\,h\frac{m_i}{\nu}h^{\epsilon}h^{-1} = CMh^{\epsilon}\frac{m_i}{\nu} < \frac{m_i}{\nu}.\label{ineq:sum_d_ij}
\end{align}
Then employing \eqref{ineq:P+} and \eqref{ineq:sum_d_ij} we get
\begin{align*}
P_i^+ & \leq \gamma_i(v_i^{\max} - v_i)\sum_{j\in \Nh^i}d_{ij}(w) 
\leq \gamma_i\frac{m_i}{\nu}(v_i^{\max} - v_i)=Q_i^+.
\end{align*}
Therefore in view of Algorithm \ref{algorithm-1} we get for $v_i>v_j$ that $\overline{\mathfrak{a}}_{ij} = R_i^+=1$.
Similar arguments may be used in case where $\overline{\mathfrak{a}}_{ij} = R_i^-=1$ and the case where $f_{ij} < 0.$
Therefore, \eqref{LP_criterion} holds and hence \eqref{eqn:linear_preserve} is satisfied. 
\end{proof}

\begin{remark}
The linearity preservation property is important in our analysis. Algorithm \ref{algorithm-1} can be implemented without defining exactly the constants $q_i$. 
In Lemma \ref{corrollary:linearity_preservation_limiters}, we consider two definitions for $q_i$ that imply the  linearity preservation property. For  $q_i$ defined by \eqref{def:q_i-2}, we assume that $h$ would be sufficiently small, such that $CM \le h^{-\epsilon}$, cf. \eqref{ineq:sum_d_ij}, where $C$ depends on the mesh and  $M$ will  on the the exact solution of the problem and will be defined later.
\end{remark}

Let us consider now the bilinear form $\widehat{d}_h(s,w;\cdot,\cdot)\,:\,{\C}\times {\C}\to{\R},$ with  $s,w\in\Sh$, defined by, for $v,z\in{\C}$, 
\begin{equation}\label{stab_term_new2}
\widehat{d}_h(s,w;v,z) := \sum_{i,j=1}^{\N}d_{ij}(w)(1 - \mathfrak{a}_{ij}(s,w))(v_i - v_j)z_i.
\end{equation} 

Then, employing $\widehat{d}_h$, we can rewrite the  AFC scheme  \eqref{ode_u_afc_2D} 
equivalently in a  variational form as: Find $u_h(t)\in \Sh$ and $c_h(t)\in \Sh$, with $u_h(0)   = u^0_h\in \Sh$, 
$c_h(0)   = c^0_h\in\Sh$, such that
\begin{equation}\begin{aligned}
(u_{h,t},\chi)_h + (\nabla u_h -\lambda  u_h\nabla c_h, \nabla \chi) + \widehat{d}_h(u_h,c_h;u_h,\chi)  & = 0,
\quad\forall\,\chi\in \Sh,\label{afc_fem_u_2D}\\
(c_{h,t},\chi)_h + (\nabla c_h,\nabla \chi) + (c_h - u_h, \chi)_h & = 0,\quad\forall\,\chi\in \Sh.
\end{aligned}\end{equation}

Note that if $\mathfrak{a}_{ij}\equiv0$, then 
$\widehat{d}_h=d_h$ and that   $\widehat{d}_h$ and  ${d}_h$ satisfy similar properties. Further, the bilinear forms $d_h$ and $\widehat d_h$
can be viewed as  $\overline{d}_h(s,w;\cdot,\cdot):{\C}\times {\C}\to{\R},$ with $s,w\in\Sh$, where
\begin{equation}\label{stab_term_afc_2D}
\overline{d}_h(s,w ;v ,z ) := \sum_{i,j=1}^{\N}\,d_{ij}(w)\rho_{ij}(s,w)(v_i - v_j)z_i,\quad\forall v,z\in{\C},
\end{equation}
with $\rho_{ij} = \rho_{ji}\in[0,1]$. Note that, for $\rho_{ij}=1$ we have $\overline{d}_h=d_h$ and for 
$\rho_{ij} = 1-\mathfrak{a}_{ij}$, we get $\overline{d}_h=\widehat d_h$.
In the sequel, we will derive various error bounds involving the bilinear form $\overline{d}_h$, which obviously will also hold for $\widehat d_h$. In view of  the symmetry of $\rho_{ij}$ and $d_{ij}$, the form 
\eqref{stab_term_afc_2D} can be rewritten as, see, e.g., \cite{gabriel2018},
\begin{equation}\label{equiv_stab_term_afc}
\overline{d}_h(s,w;v,z)  = \sum_{i<j}d_{ij}(w)\rho_{ij}(s,w)(v_i - v_j)(z_i - z_j).
\end{equation}
Hence, the schemes \eqref{low_fem_u_2D} and \eqref{afc_fem_u_2D}  
can be viewed as the following variational problem: Find $u_h(t)\in \Sh$ and $c_h(t)\in \Sh$ with $u_h(0)   = u^0_h\in\Sh$ and $c_h(0)   = c^0_h\in\Sh$, such that, 
\begin{equation}
\begin{aligned}
(u_{h,t},\chi)_h + (\nabla u_h - \lambda  u_h\nabla c_h, \nabla \chi) + \overline{d}_h(u_h,c_h;u_h,\chi)  & = 0,
\quad\forall\,\chi\in \Sh,\label{gen_fem_u_2D}\\
(c_{h,t},\chi)_h + (\nabla c_h,\nabla \chi) + (c_h - u_h, \chi)_h & = 0,\quad\forall\,\chi\in \Sh.
\end{aligned}
\end{equation}

\section{Fully discrete scheme}\label{section:fully_discrete}

 Let $\NT\in\mathbb{N}$, $\NT\ge1$, $k=T/\NT$ and $t^n=nk$, $n=0,\dots, \NT$. 
 Discretizing in time \eqref{gen_fem_u_2D} 
 with the backward Euler method we approximate 
$(u^n,c^n)=(u(\cdot,t^n),c(\cdot,t^n))$ by $(U^n,C^n)\in \Sh\times\Sh$, for $n=0,1,\dots,\NT$, such that, 
\begin{equation}
\begin{aligned}
(\overline{\partial} U^{n},\chi)_h + ( \nabla U^{n}-\lambda  U^{n}\nabla C^{n}, \nabla \chi) + 
\overline{d}_h(U^n,C^n;U^n, \chi) & = 0,\quad \forall\chi\in\Sh,\\
(\overline{\partial} C^{n},\chi)_h + (\nabla C^n,\nabla \chi) + (C^n -U^n, \chi)_h & = 0,\quad \forall\chi\in\Sh,\label{gen_fl_u_2D} 
\end{aligned}
\end{equation} 
 with $U^0   = u_h^0\in\Sh$, $C^0   = c_h^0\in\Sh$ and  $\overline{\partial} U^n = (U^n - U^{n-1})/k$.
 
Given $M>0$, let $\mathcal{B}_M\subset\Sh\times \Sh$ be such that
\begin{align}\label{bounded_S_h}
\mathcal{B}_M :=   \{(\chi, \psi)\in \Sh\times \Sh\,:\,\|\chi\|_{L^{\infty}} + \|\nabla \psi\|_{L^{\infty}} \leq M \}.
\end{align}
We will henceforth assume that $M$ is sufficiently large, specifically defined as $M = \max(2M_0, 1)$. Here, $M_0 > 0$ is independent of $h$ and is determined based on sufficiently smooth solutions $u,c$  of \eqref{Minimal_model_uc}. Consequently, for these sufficiently smooth solutions, there exists a value of $M_0 > 0$ 
such that, for $t\in[0,T]$,
\begin{equation}\label{M0-bound}
\|u(t)\|_{L^{\infty}} +\|\nabla c(t)\|_{L^{\infty}}+\|R_hu(t)\|_{L^{\infty}} + 
\|\nabla R_hc(t)\|_{L^{\infty}}\le M_0.
\end{equation}

Further, we will henceforth assume that the constants $q_i$, $i\in\Nh$, in Algorithm \ref{algorithm-1} are such that Lemma \ref{corrollary:linearity_preservation_limiters} holds. Thus if  $(U^n,C^n) \in \mathcal{B}_M$, in view of  \eqref{est_dij_2D_L_infty}, we have
\begin{equation}\label{q_max-1}
\|\q\|_{\max}\le ChM.
\end{equation}

\begin{remark}\label{q-rate}
Note that in view of \eqref{q_max-1} there exist $h_M>0$ such that for $h<h_M$, $$\Oh(h^{3/2}|\log h|+h\|\q\|_{\max})=\Oh(h^{3/2}|\log h|).$$
\end{remark}

Using a contraction argument we will show that for $(U^0,C^0)=(R_hu^0,R_hc^0)\in\mathcal{B}_M$ there exists a unique solution $(U^n,C^n)$ of \eqref{gen_fl_u_2D}, 
for $n=0,\dots,\NT$, satisfying the a priori error estimates given in the following Theorem \ref{theorem:error_estimates_low_order_2D_main}, which we  prove  at the end of this section. For the proof of  Theorem \ref{theorem:error_estimates_low_order_2D_main} and other subsequent results, we will employ the following Lemmas \ref{lemma:estimate_stab_term}-\ref{lemma:chemotaxis_bounds_2D-1}, which we show  in Section \ref{subsec:auxiliary}.

\begin{theorem}\label{theorem:error_estimates_low_order_2D_main} Let $(u,c)$ be a unique,  sufficiently smooth, solution of 
\eqref{Minimal_model_uc}, with $u\in W^{1}_{\infty} ,\,c\in W^{2}_{\infty}$. Then for sufficiently large $M>0$ and $k$, $h$ sufficiently small, $k=\Oh(h^{1+\epsilon})$ with $0<\epsilon<1$, there exists $C=C(M)>0$, independent of $k, h$ and $h_M>0$, such that for $h<h_M$, and $U^0=R_hu^0$, $C^0=R_hc^0$, we have $(U^n,C^n) \in \mathcal{B}_M$, $n\le\NT$, such that
 \begin{align*}
\|U^n-u^n\|  \le C ( k + k^{-1/2}h^2 + h^{3/2}|\log h|)\ \text{ and }\
\|C^n-c^n\|_{1} \le C ( k + k^{-1/2}h^2  +h).
\end{align*}
\end{theorem}

\begin{lemma}\label{lemma:estimate_stab_term}
There are exists a constant $C,$ independent of $h,$ such that, for $w,\widetilde  w,\psi,\chi\in \Sh$,
\begin{align*}
|d_h(w;\psi, \chi) - {d}_h(\widetilde {w};\psi,\chi)\vert \leq Ch\|\nabla (w - \widetilde {w})\|_{L^{\infty}}\|\nabla \psi\|\|\nabla \chi\|.
\end{align*}
\end{lemma}
\begin{lemma}\label{lemma:estimate_1_d_h_afc-ver2}
There are exists a constant $C,$ independent of $h,$ such that, for $w,\widehat w,\psi,\chi\in \Sh$,
\begin{align*}\
|\overline{d}_h(\psi,w;\psi, \chi) - \overline{d}_h(\psi,\widehat {w};\psi,\chi)\vert \leq Ch\|\nabla (w - \widehat {w})\|_{L^{\infty}}
\|\nabla \psi\|\|\nabla \chi\|.
\end{align*}
\end{lemma}
\begin{lemma}\label{lemma:stability_stab_term2}
Let $w\in\Sh$ and $q_i$, $i\in\Nh$, the coefficients in  Algorithm \ref{algorithm-1} such that 
Lemma \ref{corrollary:linearity_preservation_limiters} holds. 
Then, for $v\in H^{2},$ and for the bilinear form $\overline d_h$ defined in \eqref{stab_term_afc_2D} with $\rho_{ij} = \mathfrak{a}_{ij}$ or $\rho_{ij} = 1-\mathfrak{a}_{ij}$, there exists a constant $C,$ independent of $h,$ such that for all $\psi, \chi \in \Sh$, 
\begin{align*}
\vert\overline{d}_h(\psi,w;\psi,\chi)\vert & \leq C(h\|\nabla w\|_{L^{\infty}} + \|\q\|_{\max})(\|\nabla (\psi - v)\|^2 + h^2\|v\|^2_2)^{1/2}\|\nabla \chi\|. 
\end{align*}
\end{lemma}
\begin{lemma}\label{lemma:estimate_2_d_h_afc}
Let the bilinear form $\overline d_h$ be defined in \eqref{stab_term_afc_2D} with $\rho_{ij} = \mathfrak{a}_{ij}$ or $\rho_{ij} = 1-\mathfrak{a}_{ij},$ then, there exists a constant $C,$ independent of $h,$
\begin{equation}
\vert\overline{d}_h(v,w;v, \chi) - \overline{d}_h(s,w;s,\chi)\vert  \leq C(h\|\nabla w\|_{L^{\infty}} + \|\q\|_{\max})\|\nabla (v - s)\|\|\nabla \chi\|,\quad \forall w,s,v,\chi\in \Sh.\nonumber
\end{equation}
\end{lemma}
\begin{lemma}\label{lemma:chemotaxis_bounds_2D-1}
Let $v, w\in \Sh$ and $u\in W^{1}_{\infty},\,c\in W^{2}_{\infty}$. Then,  there exists a constant $C$ independent of $h$ such that for $\chi\in \Sh$,
\begin{equation*}
\vert (u\nabla c- v\nabla w,\nabla \chi)\vert  \leq C( h^{3/2}|\log h| +\|\nabla w\|_{L^{\infty}}(h^2+\|v - R_hu\|) +\|\nabla(w - R_hc)\|)\|\nabla \chi\|.
\end{equation*} 
\end{lemma}
Also for the inner product $(\cdot,\cdot)_h$ introduced in \eqref{quadrature}, the following holds.
\begin{lemma}{\normalfont{\cite[Lemma 2.3]{chatzipantelidis2012}}}\label{lemma:mass_lump_error}
Let $\varepsilon_{h}(\chi ,\psi): = (\chi , \psi) - (\chi , \psi)_h$. Then,
\begin{equation}
\vert \varepsilon_{h}(\chi ,\psi)\vert\leq Ch^{i+j}\|\nabla^i\chi \|\|\nabla^j\psi\|,\;\;\;\forall\,\chi ,\psi\in \Sh,\;\;\;\text{and}\;\;\;i,j=0,1,\nonumber
\end{equation}
where the constant $C$ is independent of $h.$
\end{lemma}
As a first step to demonstrate the a priori error estimate for $(U^n,C^n) \in \mathcal{B}_M$ we show the following theorem.
\begin{theorem}\label{theorem:error_estimates_low_order_2D-pre} 
Let $(u,c)$ be a unique, sufficiently smooth, solution of \eqref{Minimal_model_uc}, with $u\in W^{1}_{\infty},\,c\in W^{2}_{\infty}$. Also, let  $1\le n_0< \NT$ such that $(U^n,C^n)\in \mathcal{B}_M$, is  the unique solution of \eqref{gen_fl_u_2D}, for $n=0,\dots, n_0$, with $U^0=R_hu^0$ and $C^0=R_hc^0$. Then for $k$, $h$ sufficiently small, there exists $C_1=C_1(M)>0$, independent of $k, h, n_0$ such that for $n=0,\dots,n_0$, we have 
\begin{equation}\label{error_estimate_low_order_u-1}
\|R_hu^n - U^n\| + \|R_hc^n - C^n\|_{1}  \le C_1 ( k + k^{-1/2}h^2+h^{3/2}|\log h|).
\end{equation}

\end{theorem}
\begin{proof}
Let
$\theta^n  = U^n - R_hu^n$, $\rho^n  = R_hu^n-u^n$, $\zeta^n = C^n - R_hc^n$, and $\xi^n  = R_hc^n-c^n$,
for $n\ge0$.
Then,   for $\theta^n$, $n=1,\dots,\NT$, we get the following  error equation, 
\begin{equation}\label{expl_euler_theta1-1}
\begin{aligned}
&(\overline{\partial}\theta^n, \chi)_h + (\nabla\theta^n, \nabla \chi) =  - (\omega^n, \chi)_h + (\delta^n, \nabla \chi) + (\rho^n, \chi)
 - \overline{d}_h(U^n,C^n; U^n, \chi) + \varepsilon_{h}(P_hu_t^n, \chi),\, \forall\chi\in \Sh,
\end{aligned}
\end{equation}
with  
\begin{equation*}
\omega^n  = (\overline{\partial}R_hu^n - P_hu_t^n)\quad\text{ and }\quad
\delta^n  = \lambda  (u^n\nabla c^n - U^n\nabla C^n). 
\end{equation*}
Using the error estimations for $R_h$ in \eqref{ritz_projection_est2_2D}, we easily obtain
\begin{equation}\label{expl_euler_omega_2D-1}
\|\omega^n\| +\|\rho^n\|  \leq C (k + h^2),\quad \text{for } n=1,\dots,\NT.
\end{equation}
Also, in view of  Lemma \ref{lemma:chemotaxis_bounds_2D-1},  and the fact that $(U^n,C^n)\in \mathcal{B}_M$,  there exists $h_M>0$, such that for $h<h_M$,
\begin{equation}\label{expl_euler_delta_2D}
|(\delta^n,\nabla\chi)|  \leq C_M (h^{3/2}|\log h|+\|\theta^n\|+\|\nabla\zeta^{n}\|)\|\nabla\chi\|, \quad \text{for } n=1,\dots,n_0,
\end{equation}
where $C_M$ is a general positive constant that depends on $M$.

In addition, the stabilization term on the right hand side of \eqref{expl_euler_theta1-1} may be written as
\begin{align*}
-\overline{d}_h(U^n,C^n; U^n, \chi) & = (-\overline{d}_h(U^n,C^n; U^n, \chi) + \overline{d}_h(R_hu^n,C^n ; R_hu^n, \chi)) - \overline{d}_h(R_hu^n,C^n; R_hu^n, \chi) = I_1 + I_2.
\end{align*}
Employing Lemmas \ref{lemma:estimate_2_d_h_afc} and \ref{lemma:stability_stab_term2}, \eqref{q_max-1}, the fact 
that $(U^n,C^n)\in \mathcal{B}_M$ and  \eqref{ritz_projection_est2_2D}, we get 
\begin{equation}\label{est_second_case}
\begin{aligned}
|I_1| &\leq C(h\|\nabla C^n\|_{L^{\infty}} + \|\q\|_{\max})\|\nabla\theta^n\|\|\nabla \chi\|
\leq CMh\|\nabla\theta^n\|\|\nabla \chi\|,\\
|I_2| & \leq  CMh(\|\nabla (R_hu^n - u^n)\|^2 + h^{2}\|u^n\|^2_{2})^{1/2}\|\nabla \chi\|\le CMh^2\|u^n\|_2\|\nabla \chi\| 
\end{aligned}
\end{equation}

Also, employing Lemma \ref{lemma:mass_lump_error} and \eqref{L2_projection_stab}, we get 
\begin{equation}\label{estimate-eh-1}
|\varepsilon_{h}(P_hu_t^n, \chi)|\le Ch^2\| \nabla P_hu_t^n\|\|\nabla \chi\|\le Ch^2\|\nabla \chi\|.
\end{equation}

Note now that in view of the symmetry of $(\cdot,\cdot)_h$, we have
\begin{equation}\label{identity-form}
(\overline{\partial}\theta^n, \theta^n)_h=\dfrac1{2k}(\|\theta^n\|_h^2-\|\theta^{n-1}\|_h^2)+\dfrac{k}2\|\overline{\partial}\theta^n\|_h^2.
\end{equation}

Hence, choosing  $\chi= \theta^n$ in \eqref{expl_euler_theta1-1} and  using \eqref{expl_euler_omega_2D-1}--\eqref{identity-form}
we get
\begin{equation}\label{theta-n:pre}
\begin{aligned}
&\dfrac1{2k}(\|\theta^n\|_h^2-\|\theta^{n-1}\|_h^2)+\dfrac{k}2\|\overline{\partial}\theta^n\|_h^2 
 + \|\nabla\theta^n\|^2 \\
&\le C(\|\omega^n\|+\|\rho^n\|) \|\theta^n\|_h + C_M (h^{3/2}|\log h|+\|\theta^n\|+\|\nabla\zeta^{n}\|) \|\nabla \theta^n\|
+CMh\|\nabla\theta^n\|^2 \\
&\quad+CMh^2\|u^n\|_2\|\nabla\theta^n\|+ Ch^2 \| \nabla\theta^n\|.
\end{aligned}
\end{equation}
In view of Remark \ref{q-rate} there exist $h_M>0$, such that for $h<h_M$ we can eliminate $\|\nabla \theta^n\|$ in the left hand side of \eqref{theta-n:pre} and obtain
\begin{equation}\label{theta-n}
\|\theta^n\|_h^2\le \|\theta^{n-1}\|_h^2+ C_{M}k(\|\nabla\zeta^{n}\|^2+\|\theta^n\|_h^2 + {E^n}),
\end{equation}
with $E^n=O(k^2+h^{3}|\log h|^2)$.

Next, in view of the second equation of \eqref{gen_fl_u_2D}, we get the following  error equation for $\zeta^n$,
\begin{align}\label{error_eq_zeta-1}
(\overline{\partial}\zeta^n,\chi)_h + (\nabla\zeta^n,\nabla \chi) + (\zeta^n, \chi)_h & = - (\widetilde {\omega}^n, \chi)_h + (\rho^n, \chi)  + \varepsilon_{h}(\widehat{\omega}^n,\chi),
\end{align}
where
\begin{equation*}
\widetilde {\omega}^n  = (\overline{\partial}R_hc^n - P_hc_t^n) + \theta^n,\quad
\widehat{\omega}^n  =  P_hc_t^n +  R_h(c^n- u^n).
\end{equation*}
Thus for sufficiently smooth $u$ and $c$ we  get
\begin{align}\label{error_bd_zeta-1}
\|\widetilde {\omega}^n\|  \le C (k + h^2 + \|\theta^n\|_h)\quad \text{and}\quad | \varepsilon_{h} (\widehat{\omega}^n,\chi)|\le C(u,c)h^2\|\nabla\chi\|.
\end{align}
Let 
\begin{equation}\label{def:tribar}
\tribar \chi\tribar^2=\|\chi\|_h^2+\|\nabla \chi\|^2.
\end{equation}
Choosing now $\chi = \overline{\partial}\zeta^{n}$ in \eqref{error_eq_zeta-1}, using similar arguments as before and 
\eqref{error_bd_zeta-1}, we get
\begin{equation}\label{nabla-zeta-n}
\tribar\zeta^{n} \tribar^2 \leq \tribar\zeta^{n-1}\tribar^2 + C k(k^2 + k^{-1}h^4+\|\theta^n\|_{h}^2).
\end{equation}
Then combining \eqref{theta-n} and \eqref{nabla-zeta-n} we get 
\begin{equation}
\|\theta^n\|_{h}^2+\tribar\zeta^{n} \tribar^2 \leq \|\theta^{n-1}\|_{h}^2+\tribar\zeta^{n-1}\tribar^2 + C_{M}k(\|\theta^n\|_{h}^2+\tribar\zeta^{n} \tribar^2+\widetilde  E^n),
\end{equation}
with $\widetilde  E^n=O(k^2+h^{3}|\log h|^{2}+k^{-1}h^4)$.
Next,  moving $\|\theta^n\|_h^2+\tribar\zeta^{n} \tribar^2$ to the left, we have for $k$ sufficiently small and 
$n=1,\dots,n_0$,
\begin{equation}\label{theta-n-new}
\|\theta^n\|_h^2+\tribar\zeta^{n} \tribar^2\le (1+C_{M}k)(\|\theta^{n-1}\|_h^2+ \tribar\zeta^{n-1} \tribar^2)
+C_Mk\widetilde  E^n.
\end{equation}
Hence, since $\|\theta^0\|_h=\tribar\zeta^{0} \tribar=0$, summing over $n,$ we further get, for $E=\max_j{\widetilde  E^j}=O(k^2+h^{3}|\log h|^{2}+k^{-1}h^4)$, the desired estimate \eqref{error_estimate_low_order_u-1}, with $C_1$ depending on $M,$ 
\begin{align*}
\|\theta^n\|^2+\|\zeta^{n}\|_1^2&\le C
(\|\theta^n\|_h^2+\tribar\zeta^{n} \tribar^2) 
\le C_MkE\sum_{\ell=0}^n(1+C_Mk)^{n-\ell+1}\le C_1^2(k^2 
+ h^{3}|\log h|^{2}+k^{-1}h^4).
\end{align*}
\end{proof}

Then if  $\{U^n,C^n\}_{n=0}^{n_0}\in \mathcal{B}_M$, $n_0\le\NT$ is  the unique solution of \eqref{gen_fl_u_2D}, in view of Theorem \ref{theorem:error_estimates_low_order_2D-pre}, we can easily derive the a priori error bounds in Theorem \ref{theorem:error_estimates_low_order_2D_main}. 

\begin{theorem}\label{theorem:error_estimates_low_order_2D} 
Let $(u,c)$ be a unique, sufficiently smooth, solution of \eqref{Minimal_model_uc}, with $u\in W^{1}_{\infty},\,c\in W^{2}_{\infty}$. Also, let  $\{U^n,C^n\}_{n=0}^{n_0}\in \mathcal{B}_M$, $n_0\ge 1$ be  the unique solution of \eqref{gen_fl_u_2D} with $U^0=R_hu^0$ and $C^0=R_hc^0$. Then for $k$, $h$ sufficiently small, there exists $C=C(M)>0$ independent of $k, h, n_0$ such that for $n=1,\dots,n_0$,
\begin{align*}
\|U^n-u^n\|  \le C ( k + k^{-1/2}h^2 + h^{3/2}|\log h|)\ \text{ and }\
\|C^n-c^n\|_{1} \le C ( k + k^{-1/2}h^2  +h).
\end{align*}
\end{theorem}
\begin{proof}
Using the error splittings $U^n-u^n=(U^n - R_hu^n)+(R_hu^n-u^n)$, $C^n-c^n=(C^n - R_hc^n)+(R_hc^n-c^n)$,
Theorem \ref{theorem:error_estimates_low_order_2D-pre} and the approximation properties of $R_h$ in \eqref{ritz_projection_est2_2D}, 
 we can easily get the desired  error estimates. 
\end{proof}

Next, we will consider the following auxiliary set $\mathcal{\widetilde {B}}_{n}$, $n\le\NT$, defined by
\begin{equation}\label{bounded_S_h_approx}
\mathcal{\widetilde {B}}_{n} : =  \{ (\chi, \psi)\in \Sh\times \Sh:   \|\chi - R_hu^n\|_h 
+ \tribar \psi - R_hc^n\tribar \leq C_2h^{1+\widetilde \epsilon}\},
\end{equation} 
where $\tribar\cdot\tribar$ is given by \eqref{def:tribar}, $\widetilde \epsilon$ such that 
 $\Oh(h^{\widetilde\epsilon})=\Oh(h^{\min\{\epsilon,({1-\epsilon})/2\}}+h^{1/2}|\log h|)\}$ and $C_2>(2C_1+1)$, with $C_1$ given by \eqref{error_estimate_low_order_u-1}.
\begin{remark}\label{Bn-bound}
Note that in view of the inverse inequalities \eqref{eq:inverse_estimate}, for sufficiently small $h$ we have that if $(\chi,\psi)\in \mathcal{\widetilde {B}}_{n}$ then  $(\chi,\psi)\in \mathcal{ {B}}_{M}$
\end{remark}

To show that there exists a unique solution $(U^n,C^n)$ of \eqref{gen_fl_u_2D}. We will consider the following iteration operator $\Gn$ and show that it is  a contraction in $\mathcal{\widetilde {B}}_{n}\cap \mathcal{B}_M$. 
Let $\Gn=(\Gn_1,\Gn_2):\Sh\times \Sh \to \Sh \times \Sh,\,(v,w) \to (\Gn_{1}v, \Gn_{2}w)$ defined by
\begin{equation}\label{existence_linear_u}
\begin{aligned}
(\Gn_{1}v - U^{n-1}, \chi)_h +  k(\nabla \Gn_{1}v - \lambda \Gn_{1}v\nabla w, \nabla \chi)  +  
kd_h(w;\Gn_{1}v, \chi) &=   k\widetilde {d}_h(w,v;v, \chi),\quad\forall\chi\in \Sh, \\
(\Gn_{2}w - C^{n-1}, \chi)_h +  k(\nabla \Gn_{2}w, \nabla \chi)  +  k(\Gn_{2}w - v, \chi)_h  &= 0,
\quad\forall\chi\in \Sh,
\end{aligned}
\end{equation}
where,  for  $w,s\in{\Sh}$,
the bilinear form $\widetilde {d}_h(w,s;\cdot,\cdot):{\Sh}\times {\Sh}\to{\R},$   is defined by $ \widetilde {d}_h(w,s;v,z) :=  {d}_h(w;v,z)-\overline{d}_h(s,w;v,z)$,  
therefore
\begin{equation}\label{stab_term_new}
\widetilde {d}_h(w,s;v,z) =  \sum_{i,j=1}^{\N}d_{ij}(w)(1-{\rho}_{ij}(s)) (v_i - v_j)z_i.
\end{equation}
In particular,  for the low-order scheme we get $\widetilde  d_h=0$, since for this case $\rho_{ij}=1$. Note that, Lemmas \ref{lemma:estimate_1_d_h_afc} and \ref{lemma:estimate_2_d_h_afc} also hold  for
$ \widetilde {d}_h$.
Obviously, if $\Gn$ has a fixed point $(v^\star,w^\star)$, then $(U^n,C^n):=(v^\star,w^\star)$ is the solution of the discrete scheme \eqref{gen_fl_u_2D}.

We can easily rewrite \eqref{existence_linear_u} 
 in matrix formulation. 
Let $\widehat{\alb}=(\widehat{\al}_1,\dots,\widehat{\al}_{\N})^T$ and 
$\widehat{\beb}=(\widehat{\be}_1,\dots,\widehat{\be}_{\N})^T$,  the coefficients, with respect to the basis of $\Sh$, of
 $v,w\in\Sh$, respectively, and  $\widetilde {\alb}, \widetilde {\beb}\in\Rn$ the corresponding vectors 
for $\Gn_1v, \Gn_2w$, respectively. 
Then \eqref{existence_linear_u} 
can be written as
\begin{equation}\label{dasdasdasda}
\Aa \widetilde {\alb}= \ba\quad\text{and}\quad \Ab \widetilde {\beb}= \bb,
\end{equation}
where
\begin{equation}\label{dasdasdasda2}
\begin{aligned}
\Aa=\M_L +  k\,(\S - \Q_{\widehat\beb} - \D_{\widehat\beb}), &\qquad \ba=\M_L\alb^{n-1} +  k\,\overline{\bff}^{n}(\widehat\alb,\widehat\beb),\\
\Ab=\M_L +  k\,(\M_L + \S), &\qquad \bb= \M_L \beb^{n-1} + {k}\M_L{\alb}^{n-1}.
\end{aligned}
\end{equation}

Note now that in view of \eqref{zero_sum}, the matrices $\S,\,\D_{\widehat\beb}$ and $\Q_{\widehat\beb}$ have zero column sum for all $k>0.$
Also, $\Aa$, $\Ab$ have non-positive off-diagonal elements and positive diagonal elements. Then every 
 column sum of $\Aa$, $\Ab$ is positive, therefore they are strictly column diagonally dominant.
Thus for $(v,w)\in \mathcal{B}_M$ there exists a unique solution $(\Gn_1v,\Gn_2w)\in \Sh\times \Sh$ of the discrete scheme \eqref{existence_linear_u}. 

In view of Theorem \ref{theorem:error_estimates_low_order_2D-pre}, recall that if  $(U^{n-1}, C^{n-1})\in \mathcal{B}_M$ 
\begin{equation}\label{theorem:bounded_R_h-assumption}
\begin{aligned}
\|U^{n-1} - R_hu^{n-1}\| +& \| C^{n-1}- R_hc^{n-1}\|_{1}
\le C_1( k + k^{-1/2}h^2 + h^{3/2}|\log h|).
\end{aligned}
\end{equation}
Next, in the following  three lemmas we show that for $(U^{n-1}, C^{n-1})\in \mathcal{B}_M$ the iteration operator $\Gn$ is a contraction in $\mathcal{\widetilde {B}}_{n}\cap \mathcal{B}_M.$

\begin{lemma}\label{lemma-initiate-iteration} 
Let $(u,c)$ be a unique,  sufficiently smooth, solution of \eqref{Minimal_model_uc}, with $u\in W^{1}_{\infty},\,c\in W^{2}_{\infty}$. Also let $(U^{n-1}, C^{n-1})\in \mathcal{B}_M,$ such that \eqref{theorem:bounded_R_h-assumption} holds. 
Then for $k = \Oh(h^{1+\epsilon})$, with $0< \epsilon<1$, there exists exists a sufficiently large constant $C_2$, independent of $k,\,h$, such that $C_2>(2C_1+1)$ and $(U^{n-1}, C^{n-1})\in \mathcal{\widetilde {B}}_{n}.$ 
\end{lemma}

\begin{proof}
Using the stability property of $R_h$ and the fact that $k = \Oh(h^{1+\epsilon})$, we obtain the desired result
\begin{align*}
\|U^{n-1} - R_hu^{n}\|_h + \tribar C^{n-1}- R_hc^{n})\tribar
&\le \|U^{n-1} - R_hu^{n-1}\|_h + \tribar C^{n-1}- R_hc^{n-1}\tribar 
 +Ck(\|R_h\overline{\partial}u^{n}\|+\|R_h\overline{\partial}c^{n}\|_1)\\
&\le CC_1( k + k^{-1/2}h^2 + h^{3/2}|\log h|)+ C(u,c)k
\le C_2h^{1+\widetilde \epsilon}.
\end{align*}
\end{proof}

\begin{lemma}\label{theorem:bounded_R_h}
Let $(u,c)$ be a unique,  sufficiently smooth, solution of  \eqref{Minimal_model_uc}, with $u\in W^{1}_{\infty},\,c\in W^{2}_{\infty}$. Also let $(U^{n-1}, C^{n-1})\in \mathcal{B}_M$ such that  \eqref{theorem:bounded_R_h-assumption} holds 
and $(v,w)\in\mathcal{\widetilde {B}}_{n}$. Then for $k = \Oh(h^{1+\epsilon})$ with $0<\epsilon<1$, we have $(\Gn_{1}v,\Gn_{2}w)\in\mathcal{\widetilde {B}}_{n}\cap \mathcal{B}_M$.
\end{lemma}
\begin{proof}
Let $p_{n}= \Gn_{1}v - R_hu^n$, $p_{n-1} = U^{n-1}- R_hu^{n-1}$,
 $z_{n} = \Gn_{2}w - R_hc^n$, and  $z_{n-1}= C^{n-1} - R_hc^{n-1}$.
In view of the inverse inequalities \eqref{eq:inverse_estimate}, for sufficiently small $h$ we have that if $(\chi,\psi)\in \mathcal{\widetilde {B}}_{n}$ then  $(\chi,\psi)\in \mathcal{ {B}}_{M}$.

Then, in view of \eqref{existence_linear_u}, $p_n$ satisfies the following error equation for $\chi\in\Sh$,
\begin{equation}\label{pn_equation}
\begin{aligned}
(\overline{\partial} p_{n}, \chi)_h & + (\nabla p_{n}, \nabla \chi)  + d_h(w; p_{n}, \chi) = - (\omega^n, \chi)_h + (\rho^n, \chi) + (\delta^{n}, \nabla \chi) - \varepsilon_{h}(P_hu_t^{n}, \chi)
 - d_h(w; R_hu^n, \chi)\\
 & + \widetilde {d}_h(R_hc^n, R_hu^n; R_hu^n, \chi)  +\{ \widetilde {d}_h(w, v; v, \chi) - \widetilde {d}_h(R_hc^n, R_hu^n; R_hu^n, \chi)\},
\end{aligned}
\end{equation}
where
\begin{equation*}
\omega^n  = (\overline{\partial}R_hu^n - P_hu_t^{n}), \quad \rho^n  = R_hu^n-u^n, \quad \delta^{n}  = \lambda  (u^n\nabla c^n-\Gn_{1}v\nabla w). 
\end{equation*}
Using the approximation properties of $R_h$, \eqref{ritz_projection_est2_2D}, we easily obtain
\begin{equation}\label{pn-est1}
\|\omega^n\|_{h} +\|  \rho^n\|\le C( k +  h^2).
\end{equation}
Next, in view of Lemma \ref{lemma:chemotaxis_bounds_2D-1} and the fact that $(v,w)\in \mathcal{{B}}_{M}$,
and Remark \ref{q-rate}, we have for $h$ sufficiently small
\begin{equation}\label{pn-est2}
|(\delta^{n},\nabla\chi)|  \leq C (h^{3/2}|\log h|+M\|p_{n}\|+\|\nabla(w-R_hc^n)\|) \|\nabla\chi\|.
\end{equation}

Further, using the definition of $\widetilde d_h,$ i.e. \eqref{stab_term_new}, we get 

\begin{align*}
- d_h(w; R_hu^n, \chi) & + \widetilde {d}_h(R_hc^n, R_hu^n; R_hu^n, \chi)\\
& \quad = (d_h(R_hc^n; R_hu^n, \chi) - d_h(w; R_hu^n, \chi)) - \overline{d}_h(R_hu^n,R_hc^n; R_hu^n, \chi) = J_1 + J_2.
\end{align*}
In view of Lemma \ref{lemma:estimate_stab_term} 
 and the fact that $(v,w)\in \mathcal{\widetilde {B}}_{n}$, we have 
\begin{align}\label{pn-est4}
|J_1| \leq C\|\nabla (w - R_hc^n)\|\|\nabla R_hu^n\|\|\nabla \chi\| \leq CC_2h^{1+\tilde{\epsilon}}\|\nabla \chi\|.
\end{align}
Also, using Lemma \ref{lemma:stability_stab_term2}, and \eqref{q_max-1} we get
\begin{equation}\label{pn-est5}
\begin{aligned}
|J_2| & \leq  C(h\|\nabla R_hc^n\|_{L^{\infty}} + \|\q\|_{\max})(\|\nabla (R_hu^n - u^n)\|^2 + h^2\|u^n\|_{2}^2)^{1/2}\|\nabla \chi\|
 \leq CMh^{2}\|\nabla \chi\|.
\end{aligned}
\end{equation}
Next, we rewrite the last term in \eqref{pn_equation}, such that
\begin{align*}
\widetilde {d}_h(w, v; v, \chi) - \widetilde {d}_h(R_hc^n, R_hu^n; R_hu^n, \chi) & = \left(\widetilde {d}_h(w, v; v, \chi) - \widetilde {d}_h(w, R_hu^n; R_hu^n, \chi) \right)\\
& + \left( \widetilde {d}_h(w, R_hu^n; R_hu^n, \chi) - \widetilde {d}_h(R_hc^n, R_hu^n; R_hu^n, \chi) \right)  = I_1 + I_2.
\end{align*}
Note that Lemmas \ref{lemma:estimate_2_d_h_afc} and \ref{lemma:stability_stab_term2}  are valid also for  $\widetilde  d_h$.
Using Lemma \ref{lemma:estimate_2_d_h_afc}, the inverse estimate, \eqref{eq:inverse_estimate}, the fact that $(v,w)\in\mathcal{\widetilde {B}}_{n}$ and \eqref{q_max-1}, we obtain
\begin{equation}\label{pn-est6}
|I_1|  \leq C(h\|\nabla w\|_{L^{\infty}} + \|\q\|_{\max})\|\nabla (v - R_hu^n)\|\|\nabla \chi\| \leq Ch^{-1}hM \|v - R_hu^n\|\|\nabla \chi\|
 \leq CMC_2h^{1+\widetilde{\epsilon}}\|\nabla \chi\|.
\end{equation}
Next,  in view of Lemma \ref{lemma:stability_stab_term2} and  \eqref{q_max-1}, we get
\begin{align}\label{pn-est7}
|I_2|  \leq CMh(\|\nabla (R_hu^n - u^n)\|^2 + h^2\|u^n\|^2_2)^{1/2}\|\nabla \chi\| 
\le CMh^{2}\|\nabla \chi\|.
\end{align}

Choosing now $\chi = p_{n}$ in \eqref{pn_equation}, employing the corresponding identity as in 
\eqref{identity-form} for $p_n$, combining the previous estimations \eqref{estimate-eh-1}, 
\eqref{pn-est1}-\eqref{pn-est7},  
and eliminating $\|\nabla p_{n}\|$ we get 
\begin{align}
\|p_{n}\|_{h}^2 &\le  \|p_{n-1}\|_{h}^2 +  C_{M}k (\|p_{n}\|_{h}^2 + k^2+h^{2+2\widetilde \epsilon})\label{p-n}.
\end{align}
with $C_{M}$ denoting a constant that depends on $C_2$ and $M$.

Next, in view of \eqref{existence_linear_u}, we have the following error equation for $z_n$,
\begin{equation}\label{qn_equation}
(\overline{\partial} z_{n}, \chi)_h + (\nabla z_{n}, \nabla \chi) + (z_{n}, \chi)_h   = - (\omega^{n}, \chi)_h + (\rho^n, \chi) - \varepsilon_{h}(\widehat{\omega}^n, \chi),
\end{equation}
where
\begin{equation*}
\omega^{n}  = u^n - v + (\overline{\partial}R_hc^n - P_hc_t^n), \quad\text{and}\quad \widehat{\omega}^n = R_h(c^n - u^n) + P_hc_t^n.
\end{equation*}
Employing now \eqref{ritz_projection_est2_2D} and the fact that $(v,w)\in \widetilde {\mathcal{B}}_n$, we obtain
\begin{equation*}
\|\omega^{n}\|_{h} + \|\rho^{n}\| \leq C(k+h^2)+C_2h^{1+\widetilde \epsilon}\quad\text{and}\quad | \varepsilon_{h} (\widehat{\omega}^n,\chi)|\le Ch^2\|\nabla\chi\|.
\end{equation*}
Choosing now $\chi = \overline{\partial}z_{n}$ in \eqref{qn_equation}
and  using similar arguments as before, we get
\begin{align}
\tribar z_{n} \tribar^2 &\leq \tribar z_{n-1}\tribar^2 + C_{M}k(k^2 +h^{2+2\widetilde \epsilon})+Ch^4.\label{q-n}
\end{align}
Then combining \eqref{p-n} and \eqref{q-n}, we have for $k$ sufficiently small
\begin{align*}
\|p_{n}\|_{h}^2  + \tribar z_{n}\tribar^2 \le (1+C_{M}k) (\|p_{n-1}\|_{h}^2  + \tribar z_{n-1}\tribar^2)+ C_{M}kE,
\end{align*}
with $E=\Oh(k^2+h^{2+2\widetilde \epsilon}+ k^{-1}h^{4})$.
Then for $k=\Oh(h^{1+\epsilon})$ sufficiently small, we obtain
\begin{align*}
\|p_{n}\|_{h}^2  + \tribar z_{n}\tribar^2 \le 
C_2^2h^{2+2\widetilde \epsilon}.
\end{align*}
Hence, $(\Gn_{1}v,\Gn_{2}w)\in \mathcal{\widetilde {B}}_{n}$.
Finally, employing the inverse inequality \eqref{eq:inverse_estimate}, the fact that 
$(\Gn_{1}v,\Gn_{2}w)\in \mathcal{\widetilde {B}}_{n}$ and \eqref{ritz_projection_stab}, we get sufficiently small $h$,
\begin{align*}
\| \Gn_{1}v \|_{L^{\infty}}+\|\nabla \Gn_{2}w \|_{L^{\infty}} & \le Ch^{-1}(\| p_{n}\|_h+\tribar z_{n}\tribar) 
+ \| R_hu^n\|_{L^{\infty}} +\|\nabla R_hc^n\|_{L^{\infty}} 
\le Ch^{\widetilde \epsilon}+M_0\le M.
\end{align*}
Therefore $(\Gn_{1}v, \Gn_{2}w)\in \mathcal{B}_{M},$ which concludes the proof.
\end{proof}

Now, we can prove the existence and uniqueness of the solution of \eqref{gen_fl_u_2D}.
 
\begin{theorem}\label{theorem:existence_fully_discrete_afc_2D} Let $(u,c)$ be a unique,  sufficiently smooth, 
solution of  \eqref{Minimal_model_uc}, with $u\in W^{1}_{\infty},\,c\in W^{2}_{\infty}$. Let the correction factors $\mathfrak{a}_{ij}^n,$ for $i,j=1,\dots,\N,$ be defined as in Algorithm \ref{algorithm-1}. If $(U^{n-1},C^{n-1}), (v,w)\in\mathcal{B}_{M},$ such that \eqref{theorem:bounded_R_h-assumption} holds. Then for $h$ sufficiently small,  $k = \Oh(h^{1+\epsilon})$ 
 with $\epsilon>0$, there exists a unique solution $(U^n, C^n)\in \mathcal{{B}}_{M}$ of the fully-discrete scheme \eqref{gen_fl_u_2D}.
\end{theorem}
\begin{proof}
Obviously, in view of Lemmas \ref{lemma-initiate-iteration} and \ref{theorem:bounded_R_h}, 
starting with $(v_0,w_0)=(U^{n-1},C^{n-1})$ through $\Gn$, we obtain a sequence of elements 
$(v_{j+1},w_{j+1})=(\Gn_1v_j,\Gn_2w_j)\in\mathcal{\widetilde {B}}_{n}\cap\mathcal{B}_{M}$, $j\ge0$. 

To show existence and uniqueness of $(U^n, C^n)\in \mathcal{\widetilde {B}}_{n}\cap\mathcal{B}_{M}$ it suffices to show that there exists $0<L<1$, such that 
\begin{equation*}
\|\Gn_1v-\Gn_1\widehat{v}\|_h+\tribar\Gn_2w-\Gn_2\widehat w\tribar\le L (\|v-\widehat{v}\|_h+\tribar w-\widehat w\tribar),\quad \forall (v,w),(\widehat v,\widehat w)\in \mathcal{\widetilde {B}}_{n}\cap\mathcal{B}_{M}.
\end{equation*}

Let $(v,w), ({\widehat v},{\widehat w})\in \mathcal{\widetilde {B}}_{n}\cap \mathcal{{B}}_{M},$ with $\widetilde {v}:= v - {\widehat v}$ and $\widetilde {w}: = w - {\widehat w}.$ In view of \eqref{existence_linear_u}, we have
\begin{equation}\label{Lip-v-1}
\begin{aligned}
&(\Gn_{1}\widetilde {v}, \chi)_h  +  k(\nabla \Gn_{1}\widetilde {v}, \nabla \chi) +k\,d_h(w;\Gn_{1}\widetilde {v}, \chi) \\
&=k \lambda (\Gn_{1}{v}\nabla  w-\Gn_{1}{\widehat v}\nabla \widehat {w}, \nabla \chi) 
-k\,\{d_h(w;\Gn_{1}\widehat v, \chi) -  d_h({\widehat w};\Gn_{1}{\widehat v}, \chi) \}
  +   k\,\{\widetilde {d}_h(w, v;v, \chi) -   \widetilde {d}_h({\widehat w}, {\widehat v};{\widehat v}, \chi)\}\\
&\quad =I_1+I_2+I_3. 
\end{aligned}
\end{equation}
We can rewrite $I_1$  and $I_3$, in the following way
\begin{align*}
I_1 & = k \lambda (\Gn_{1}{v}\nabla \widetilde w, \nabla \chi) 
+k \lambda (\Gn_{1}{\widetilde v}\nabla \widehat w, \nabla \chi) 
 =I_1^1+I_1^2,\\
I_3 &= k\{\widetilde {d}_h({w}, v; v, \chi) - \widetilde {d}_h({\widehat w}, v; v, \chi)\}
 + k\{\widetilde{d}_h({ \widehat w}, v; v, \chi) - \widetilde {d}_h({\widehat w}, {\widehat v}; {\widehat v}, \chi)\}
 =I_3^1+I_3^2. 
\end{align*}
Then, in view of the fact that 
$(\Gn_1{\widehat v},\Gn_2{\widehat w})\in\mathcal{\widetilde {B}}_{n}\cap\mathcal{B}_{M}$ and \eqref{eq:inverse_estimate}
\begin{equation}\label{Gn-bound}
\|\nabla \Gn_{1}{\widehat v}\|  \leq C\|\nabla (\Gn_{1}{\widehat v} - R_hu^n)\|  + \|\nabla R_hu^{n}\| 
 \leq Ch^{-1}\|\Gn_{1}{\widehat v} - R_hu^n\|_h  + \|\nabla R_hu^{n}\|
 \leq Ch^{\tilde{\epsilon}} + M_0 \leq M.
\end{equation}
Similarly, we obtain $\|\nabla v\| \leq M.$ Using the fact that $(v,w)$, $({\widehat v},{\widehat w})$, $(\Gn_1v,\Gn_2w)\in \mathcal{\widetilde {B}}_{n}\cap\mathcal{B}_{M} $ and Lemmas 
\ref{lemma:estimate_stab_term}, \ref{lemma:estimate_2_d_h_afc} and \ref{lemma:estimate_1_d_h_afc-ver2}, \eqref{q_max-1}, \eqref{Gn-bound} and \eqref{eq:inverse_estimate},  we  get
\begin{equation}\label{I1-estimations}
\begin{aligned}
|I_1| & \le CMk(\|\nabla \widetilde {w}\| + \|\Gn_{1}\widetilde {v}\|)\|\nabla \chi\|,\quad  |I_3|  \leq CMk
(\|\nabla \widetilde {w}\|+\|\widetilde {v}\|)\|\nabla \chi\|,\\
|I_2| & \le Ck\|\nabla \widetilde {w}\|\|\nabla \Gn_{1}{\widehat v}\|\|\nabla \chi\|\le  CMk\|\nabla \widetilde {w}\|\|\nabla \chi\|.\\
\end{aligned}
\end{equation}
 Choosing $\chi = \Gn_{1}\widetilde {v}$ in \eqref{Lip-v-1}, using \eqref{I1-estimations} and eliminating $\|\nabla \Gn_{1}\widetilde {v}\|$, we get
\begin{equation*}
\|\Gn_{1}\widetilde {v}\|_h^2   \leq C_{M} k( \|\Gn_{1} \widetilde {v}\|_h^2 + \|\nabla \widetilde {w}\|^2  +  \|\widetilde {v}\|_h^2). 
\end{equation*}
Then for sufficiently small $k$, we obtain
\begin{equation}\label{exi_eq1}
\|\Gn_{1}\widetilde {v}\|_h^2   \leq C_{M} k(  \|\nabla \widetilde {w}\|^2  +  \|\widetilde {v}\|_h^2).
\end{equation}

Next, in view of \eqref{existence_linear_u}, we get
\begin{equation}\label{Lip-w-1}
\begin{aligned}
(\Gn_{2}\widetilde {w}, \chi)_h  +  k\,(\nabla \Gn_{2}\widetilde {w}, \nabla \chi) +k(\Gn_{2}\widetilde {w}, \chi)=k(\widetilde  v,\chi)
\end{aligned}
\end{equation}
Choosing $\chi= \Gn_{2}\widetilde {w}$ in \eqref{Lip-w-1} we get
\begin{align*}
\|\Gn_{2}\widetilde {w}\|^2 +   k\,\tribar\Gn_{2}\widetilde {w}\tribar^2   & \leq  Ck\,\|\widetilde {v}\|\|\Gn_{2}\widetilde {w}\|. 
\end{align*}
which after eliminating $\|\Gn_{2}\widetilde {w}\|$ gives
\begin{equation}\label{exi_eq2}
\tribar\Gn_{2}\widetilde {w}\tribar^2  \leq  Ck\|\widetilde {v}\|_h^2. 
\end{equation}

Thus combining \eqref{exi_eq1} and \eqref{exi_eq2}, we get
\begin{equation*}
\|\Gn_{1}\widetilde {v}\|_h^2 +  \tribar\Gn_{2}\widetilde {w}\tribar^2   \leq C k( \| \widetilde {v}\|_h^2
+\tribar \widetilde {w}\tribar^2). 
\end{equation*}
Therefore, for $k$ sufficiently small such that $C k < 1$, the sequence $(v_{j},w_{j})\to (U^n,C^n)\in\Bn\cap\BM$, 
$j\to\infty$ and $(U^n,C^n)$ is the unique solution  of the fully-discrete scheme \eqref{gen_fl_u_2D}.
\end{proof}

Finally, combining the previous results we can show Theorem \ref{theorem:error_estimates_low_order_2D_main}.
\begin{proof}[Proof of Theorem \ref{theorem:error_estimates_low_order_2D_main}] 
Obviously, in view of Lemma \ref{lemma-initiate-iteration}, for $U^0=R_hu^0$, $C^0=R_hc^0$, $(U^0,C^0) \in \mathcal{B}_M\cap \mathcal{\widetilde {B}}_{n}$.
Employing now Theorem \ref{theorem:existence_fully_discrete_afc_2D}, we get that there exists a unique solution $(U^n,C^n)$ of \eqref{gen_fl_u_2D}, for $n=0,\dots,\NT$. Then, using Theorem \ref{theorem:error_estimates_low_order_2D} we obtain the desired error bounds.
\end{proof}

\section{Positivity}\label{section:positivity}
In this section we will demonstrate that the solution $(U^n,C^n)$ of the 
the fully discrete scheme \eqref{gen_fl_u_2D} is positive if for the initial approximations $(U^{0},C^{0})$ are positive. A similar assumption for the strict positivity of the initial approximations also assumed in \cite{filbet2006,li-shu2017}.

The fully discrete scheme \eqref{gen_fl_u_2D} may be expressed by splitting the bilinear form $\overline d_h$ in a similar way as the iteration scheme \eqref{existence_linear_u},
\begin{align*}
(U^n - U^{n-1}, \chi)_h +  & k\,(\nabla U^n, \nabla \chi)  - \lambda\, k\,(U^n\nabla C^n, \nabla \chi)  +  k\,d_h(C^n;U^n, \chi) = k\,\widetilde {d}_h(C^n,U^n;U^n, \chi),\quad\forall\chi\in\Sh, \\
(C^n - C^{n-1}, \chi)_h +  &k\,(\nabla C^n, \nabla \chi)  +  k(C^n - U^n, \chi)_h  = 0,\quad\forall\chi\in\Sh.
\end{align*} 
We can easily rewrite \eqref{gen_fl_u_2D} in matrix formulation. For this, we introduce 
the following notation. Let ${\alb}^n=({\al}^n_1,\dots,{\al}^n_{\N})^T$ and ${\beb}^n=({\be}^n_1,\dots,{\be}^n_{\N})^T$,  the coefficients, with respect to the basis of $\Sh$, of $U^n,C^n\in\Sh$, respectively.
Then \eqref{gen_fl_u_2D} can be written as
\begin{equation}\label{dasdasdasda-1}
\Aan {\alb}^n= \ban\quad\text{and}\quad \Abn {\beb}^n= \bb,
\end{equation}
where
\begin{equation}\label{dasdasdasda2-1}
\begin{aligned}
\Aan=\M_L +  k\,(\S  - \Q_{\beb^n} - \D_{\beb^n}), &\ \ban=\M_L\alb^{n-1} +  k\,\overline{\bff}^{n}(\alb^n,\beb^n),\\
\Ab=\M_L +  k\,(\M_L + \S ), &\quad \bb= \M_L \beb^{n-1} + {k}\M_L{\alb}^{n}.
\end{aligned}
\end{equation}

Next, we will assume that $U^0,\,C^0>0$ and prove that $U^n,\,C^n>0,\,n=1,\ldots,\NT,$ for the solutions of the discrete scheme \eqref{gen_fl_u_2D} are non-negative,  for small $k$ and $k = \Oh(h^{1+\epsilon})$ for $0<\epsilon<1$.

\begin{theorem}\label{theorem:positivity_afc_2D} 
Let the correction factors $\mathfrak{a}_{ij}^n,$ for $i,j=1,\dots,\N$, be defined as in Algorithm \ref{algorithm-1}. Then for positive initial approximation $(U^0,C^0)$, the solution of the discrete scheme \eqref{gen_fl_u_2D} $\{(U^n, C^n)\}_{n=0}^{\NT},$ is also positive, for small $k,\,h,$ and $k = \Oh(h^{1+\epsilon})$ for $0<\epsilon<1.$
\end{theorem}
\begin{proof}
Note that in order to show the positivity of $\{(U^n,C^n)\}_{n=1}^{\NT}$ it suffices to show that the vectors $\alb^{n}$, $\beb^{n}$ of the coefficients of $\{U^n, C^n\}_{n=0}^{\NT}$ with respect to the basis of $\Sh$,  
are positive $\alb^{n} > 0,\,\beb^{n}\geq 0,$ element-wise. We will show this by induction. 
Since $U^0>0,\,C^0>0$, it suffices show that $U^n>0,\,C^n>0$ if $U^{n-1}>0,\,C^{n-1}>0$.
The assumption that $U^{n-1}>0,\,C^{n-1}>0$ implies that $\alb^{n-1}>0,\,\beb^{n-1}>0$, respectively.
Therefore, in order to show the desired result we will show that $\alb^{n}>0,\,\beb^{n}>0$.
 
The fully discrete scheme \eqref{gen_fl_u_2D} can equivalently be written in matrix formulation, i.e., \eqref{dasdasdasda-1}--\eqref{dasdasdasda2-1}. From the first system, i.e., $\Aan {\alb}^n= \ban$, we obtain,  for  $i\in \Nh$,
\begin{align*}
(m_i + &k(s_{ii} - \tau_{ii}^n - d_{ii}^n))\al_i^n +  k\sum_{j\in \Nh^i}( s_{ij} - \tau_{ij}^n - d_{ij}^n)\al_j^n 
= m_i\al_i^{n-1} +k\sum_{j\in \Nh^i}\mathfrak{a}_{ij}^nd_{ij}^n(\al_i^n - \al_j^n),
\end{align*}
where $\Q_{\beb^n}=(\tau_{ij}^n)_{i,j=1}^{\N} $ and $\D_{\beb^n}=(d_{ij}^n)_{i,j=1}^{\N} $. If $\al_i^n>0$, then we have the desired result. Therefore, let us assume that $\al_i^n \leq 0$ and $\al_i^n = \min_{j\in\Nh^i}\al_j^n$. Using the zero row sum property of $\D_{\beb^n}$, see also \eqref{D_def}, we obtain
\begin{align}\label{pos_1}
(m_i + &k(s_{ii} - \tau_{ii}^n ))\al_i^n +  k\sum_{j\in \Nh^i}( s_{ij} - \tau_{ij}^n )\al_j^n 
= m_i\al_i^{n-1} - k\sum_{j\in \Nh^i}(1 - \mathfrak{a}_{ij}^n)d_{ij}^n(\al_i^n - \al_j^n) \geq m_i\al_i^{n-1},
\end{align}
where the latter holds since $d_{ij}^n \geq 0,\,j\neq i,\,\mathfrak{a}_{ij}^n\in [0,1]$ and $\al_i^n < \al_j^n$ for all $j\in\Nh^i.$ Then, for small $k$ and $k=\mathcal{O}(h^{1+\epsilon}),\,0<\epsilon<1,$ in view of Theorem \ref{theorem:existence_fully_discrete_afc_2D}, yields that $\|\nabla C^n\|_{L^{\infty}} \leq M,$ with $M>0$ for all $n=0,\ldots,\N.$ Using the latter a priori estimate together with \eqref{D_def}, \eqref{est_dij_2D_L_infty}, we obtain for $j\neq i,$ that $s_{ij} - \tau_{ij}^n \leq s_{ij} + d_{ij}^n \leq s_{ij} + C_M\,h \leq 0$ for $h < \frac{|s_{ij}|}{C_M}.$ Note that, in two dimensions the elements of stiffness matrix $\S$ are independent of $h.$ Thus,  since $\al_i^n < \al_j^n,\,j\in\Nh^i,$ and for $h$ such that $s_{ij} - \tau_{ij}^n\leq 0,\,j\neq i,$ yields
\begin{align*}
(m_i + k(s_{ii} - \tau_{ii}^n ))\al_i^n +  k\sum_{j\in \Nh^i}( s_{ij} - \tau_{ij}^n )\al_i^n \geq (m_i + k(s_{ii} - \tau_{ii}^n ))\al_i^n +  k\sum_{j\in \Nh^i}( s_{ij} - \tau_{ij}^n )\al_j^n \geq m_i\al_i^{n-1},
\end{align*}
where the last estimate is due to \eqref{pos_1}.
Using the zero row sum property of $\S$, we further obtain,
\begin{align*}
\left( m_i - k\tau_{ii}^n -  k\sum_{j\in \Nh^i}\tau_{ij}^n\right) \al_i^n \geq m_i\al_i^{n-1}> 0.
\end{align*}
Therefore, $(m_i-  k\,\sum_{j\in \Nh^i\cup\{i\}}\tau_{ij}^n)\al_i^n>0$, since $\al_i^{n-1}>0.$ On the other hand, for  $k=O(h^{1+\epsilon}),\,0<\epsilon<1$, for sufficiently small $h$, there exists constants $c_1,c_2>0$ such that
\begin{align*}
m_i-  k\,\sum_{j\in \Nh^i\cup\{i\}}\tau_{ij}^n & \geq c_1h^2 - c_2Mh^{2+\epsilon} > 0.
\end{align*}
Hence, we arrive in a contradiction, thus $\al_i^n>0,$ for $i\in\Nh.$ Note that, prove the positivity of the term inside the parenthesis, we have used also the uniform boundedness of $\|\nabla C^n\|_{L^\infty} \leq M,\,M>0,$ for all $n=0,\ldots,\N.$

To show the positivity of $\beb^n$ given that $\beb^{n-1}>0,$ if suffices in view of \eqref{dasdasdasda2-1} to ensure that $\Abn$ has positive inverse and the vector $\bb$ is also positive. In particular, the matrix $\Abn$ is M-matrix for all $k$, since its off-diagonal elements are non-positive and the diagonal elements are positive. Further, it can be proved that is invertible, since it is strict diagonally dominant by columns. The positivity of  $\bb$ follows for the positivity of $\beb^{n-1},\,\alb^n.$
\end{proof}

Note that the latter theorem is also true when $C^0\geq 0.$ For the discrete approximation $U^n$ of $u^n$ we can show the following conservation property.

\begin{lemma}\label{lemma:afc_order_fully_discrete_total_mass_2D}
Let $\{U^n,\,C^n\}_{n=0}^{\NT}$ be a solution of the fully-discrete scheme \eqref{gen_fl_u_2D} 
Then, we have
\begin{equation}\label{afc_conservation_total_mass_2D}
(U^n,1)_h = (U^0,1)_h,\;\;\;\;\;\;\;\text{for all}\;\;n= 1,\dots,\NT.
\end{equation}
In addition, is conserved in $L^{1}$--norm for small $k$ and $k = \Oh(h^{1+\epsilon})$ for $0 <\epsilon < 1$ and  $\|\q\|_{\max} = \Oh(h)$ i.e.,
\begin{equation}\label{afc_conservation_total_massL1_2D}
\|U^n\|_{L^{1}} = \|U^0\|_{L^{1}},\;\;\text{for all}\;\;n = 1,\dots, \NT.
\end{equation}
\end{lemma}
\begin{proof}

We can easily see that choosing  $\chi=1\in\Sh$ in \eqref{gen_fl_u_2D} we get that
\begin{equation}
(\overline{\partial} U^{n},1)_h =0.
\end{equation} 
Thus we can easily obtain the conservation property \eqref{afc_conservation_total_mass_2D}.
Then, combining \eqref{afc_conservation_total_mass_2D} and positivity of $U^n>0$, $n=0,\dots,\NT$  for  $U_0>0$, we get the conservation in $L^{1}$--norm.
\end{proof}

\section{Auxiliary results}\label{subsec:auxiliary}
Here in this section we will prove some auxilliary results and Lemmas \ref{lemma:estimate_stab_term}-\ref{lemma:chemotaxis_bounds_2D-1}.
Using the following lemma  we have that the bilinear form $d_h$, introduced in \eqref{stab_term}, and hence also $\widehat{d}_h$, defined in \eqref{stab_term_new2}, induces a seminorm on $\C$.
\begin{lemma}{\normalfont{\cite[Lemma 3.1]{gabriel2016}}}\label{lemma:seminorm}
Consider any $\varpi_{ij} = \varpi_{ji}\geq 0$ for $i,j=1,\dots,\N.$ Then,
\begin{equation*}
\sum_{i,j=1}^{\N}v_i\varpi_{ij}(v_i - v_j) = \sum_{\substack{i,j=1,\, i<j}}^{\N}\varpi_{ij} (v_i-v_j)^2\geq 0,\quad \forall v_1,\dots,v_{\N}\in\mathbb{R}.
\end{equation*}
\end{lemma}
Therefore, $\overline{d}_h(s,w;\cdot,\cdot):{\C}\times {\C}\to{\R},$ with $s,w\in\Sh$, is a non-negative symmetric bilinear form 
which satisfies the Cauchy-Schwartz's inequality,
\begin{equation} \label{Schwartz_ineq_afc_2D}
|\overline{d}_h(s,w;v,z)|^2\leq \overline{d}_h(s,w;v,v)\, \overline{d}_h(s,w;z,z),\quad\forall v,z\in{\C},
\end{equation}  
and thus induces a seminorm on ${\C}$.

\begin{lemma}\label{lemma:stability_stab_term}
There are exists a constant $C,$ independent of $h,$ such that for all $w,s,\psi,\chi\in \Sh,$
\begin{align*}
\vert\overline{d}_h(s,w;\psi,\chi)\vert & \leq Ch\|\nabla w\|_{L^{\infty}}\|\nabla \psi\|\|\nabla \chi\|.
\end{align*}
\end{lemma}
\begin{proof}
Recall that \eqref{est_dij_2D_L_infty} gives 
\begin{equation*}
|d_{ij}(w)|  \leq Ch\|\nabla w\|_{L^{\infty}}, \quad\forall i,j\in\Nh.
\end{equation*}

Then, in view of \eqref{equiv_stab_term_afc} and the fact that $\rho_{ij}\in[0,1]$, we get
\begin{align}
|\overline{d}_h(s,w;\chi,\chi)| &\le  Ch\|\nabla w\|_{L^{\infty}}\sum_{K\in\Th}\sum_{i,j\in \Nh(K)} ( \chi_i - \chi_j)^2
 \le Ch\|\nabla w\|_{L^{\infty}}\sum_{K\in\Th}\|\nabla \chi\|_{L^{2}(K)}^2\notag\\
&\le Ch\|\nabla w\|_{L^{\infty}}\|\nabla \chi\|^2.\label{d_h_est_2D-2}
\end{align}
Therefore, using \eqref{Schwartz_ineq_afc_2D} we get the desired result.
\end{proof}

Next, we show two  Lipschitz-like estimates for $\overline{d}_h$.

\begin{lemma}\label{lemma:estimate_1_d_h_afc}
There are exists a constant $C,$ independent of $h,$ such that, for $w,\widetilde  w,\psi,\chi\in \Sh$,
\begin{align*}\
|\overline{d}_h(\psi,w;\psi, \chi) - \overline{d}_h(\psi,\widetilde {w};\psi,\chi)\vert \leq Ch(\|\nabla (w - \widetilde {w})\|_{L^{\infty}}+\|\nabla w\|_{L^{\infty}})
\|\nabla \psi\|\|\nabla \chi\|.
\end{align*}
\end{lemma}

\begin{proof} In view of \eqref{equiv_stab_term_afc}, we have
\begin{align*}
& \overline{d}_h(\psi,w; \psi, \chi) - \overline{d}_h(\psi,\widetilde {w}; \psi, \chi) 
= \sum_{i<j}(d_{ij}(w)\rho_{ij}(\psi,w) - d_{ij}(\widetilde {w})\rho_{ij}(\psi,\widetilde{w}))(\psi_i - \psi_j)(\chi_i-\chi_j).
\end{align*}
Also, using \eqref{D_def}, we get
\begin{align*}
&d_{ij}(\widetilde {w})\rho_{ij}(\psi,\widetilde{w}) - d_{ij}(w)\rho_{ij}(\psi,w) 
= (d_{ij}(\widetilde {w})- d_{ij}(w))\rho_{ij}(\psi,\widetilde{w}) + d_{ij}(w)(\rho_{ij}(\psi,\widetilde{w}) - \rho_{ij}(\psi,w))\\ 
&\qquad =  \{\max\lbrace - \tau_{ij}(\widetilde {w}),0,- \tau_{ji}(\widetilde {w}) \rbrace - \max\lbrace - \tau_{ij}(w),0,- \tau_{ji}(w) \rbrace\} \rho_{ij}(\psi,\widetilde{w}) 
+ d_{ij}(w)(\rho_{ij}(\psi,\widetilde{w}) - \rho_{ij}(\psi,w)) \\
&\qquad \leq \vert\tau_{ij}(w) - \tau_{ij}(\widetilde {w})\vert + \vert\tau_{ji}(w) - \tau_{ji}(\widetilde {w})\vert +2d_{ij}(w).							
\end{align*}
In addition, similar to \eqref{est_dij_2D_L_infty}, there exists a constant $ C_\gamma>0$ independent of $h$, such that
\begin{equation}\begin{aligned}\label{diff-tau_ij}
|\tau_{ij}(w) - \tau_{ij}(\widetilde {w})| &\le \sum_{K\in\omega_i}
\int_{K} |(\nabla (w - \widetilde {w})\cdot \nabla \phi_i)\phi_j|\,dx
 \le C\|\nabla (w - \widetilde {w})\|_{L^{\infty}(\omega_i)}\sum_{K\in\omega_i}h_K\\
 &\le C_\gamma h\|\nabla (w - \widetilde {w})\|_{L^{\infty}(\omega_i)}.
\end{aligned} \end{equation}

Therefore, in view of \eqref{diff-tau_ij} and \eqref{est_dij_2D_L_infty}, we obtain
\begin{align*}
\vert \overline{d}_h(\psi,w;\psi, \chi) - \overline{d}_h(\psi,\widetilde {w};\psi,\chi)\vert & 
\leq Ch(\|\nabla (w - \widetilde {w})\|_{L^{\infty}} + \|\nabla w\|_{L^{\infty}})\sum_{i<j}\vert \psi_i - \psi_j\vert\vert\chi_i - \chi_j\vert\\
 &\leq Ch(\|\nabla (w - \widetilde {w})\|_{L^{\infty}}+\|\nabla w\|_{L^{\infty}})\|\nabla \psi\|\|\nabla \chi\|,
\end{align*}
which gives the desired bound.
\end{proof}

Next, similarly as in Lemma \ref{lemma:estimate_1_d_h_afc} we get the following estimate for the stabilization term $d_h$ defined in \eqref{stab_term}.

Next, we show that the correction factor functions 
$\overline{\mathfrak{a}}_{ij}$ can be written similarly as in  {\cite[Lemma 4.1]{gabriel2016}}. 
\begin{lemma}\label{lemma:limiters_estimate}
Let the triangulation $\Th$ satisfy Assumption \ref{mesh-assumption}. Then the correction factor functions 
$\overline{\mathfrak{a}}_{ij}$,  defined by Algorithm \ref{algorithm-1}, can be written as 
\begin{equation}\label{correction_factors_formulation}
\overline{\mathfrak{a}}_{ij}(\alb,\beb) = \frac{A_{ij}(\alb,\beb)}{\vert\al_j - \al_i\vert +B_{ij}(\alb,\beb)}, \quad i,j\in\Nh, 
\end{equation}
with  $\alb=(\al_1,\dots,\al_{\N})^T,\,\beb=(\be_1,\dots,\be_{\N})^T\in \mathbb{R}^{\N}$, $\al_j\neq \al_i$, and $A_{ij}$ 
and $B_{ij}$ non-negative functions which are continuous functions  in $\alb$, with $\al_i\neq\al_j$. Further 
$A_{ij}$ 
and $B_{ij}$ are Lipschitz-continuous with respect to $\alb$  in the sets $\{\alb\in\Rn:\al_i>\al_j\}$
 and $\{\alb\in\Rn:\al_i<\al_j\}$, 
 \begin{align*}
 |A_{ij}(\alb,\beb)-A_{ij}(\widetilde\alb,\beb)|&\le \Lambda_{ij}^A(\beb,q_i)\sum_{\ell\in\Nh(\omega_i)}|\al_\ell-\widetilde\al_\ell|\\
\quad \text{ and }\quad
  |B_{ij}(\alb,\beb)-B_{ij}(\widetilde\alb,\beb)|&\le \Lambda_{ij}^B(\beb,q_i)\sum_{\ell\in\Nh(\omega_i)}|\al_\ell-\widetilde\al_\ell|.
 \end{align*}
  where  
  \begin{equation*}
  \Lambda_{ij}^A(\beb,q_i)= \Lambda_{ij}^B(\beb,q_i)+\frac{q_i}{d_{ij}(\beb)},\quad
  \Lambda_{ij}^B(\beb,q_i)=\frac{\max_{1\leq j\leq \N}d_{ij}(\beb) }{d_{ij}(\beb)}.
  \end{equation*}
 \end{lemma}
 
 \begin{proof}
 Following {\normalfont{\cite[Lemma 4.1]{gabriel2016}}} we define  $A_{ij}$ and $B_{ij}$, for $i,j\in\Nh$ and $d_{ij} > 0,$ as
\begin{align}\label{def_A_ij}
A_{ij}(\alb,\beb) := \frac{1}{d_{ij}(\beb)}\begin{cases}
\min\{ -P_{i}^-(\alb,\beb), -Q_{i}^-(\alb)\}, & \quad\text{ if }\ \al_i < \al_j,\\
\min\{ P_{i}^+(\alb,\beb), Q_{i}^+(\alb)\}, & \quad\text{ if }\ \al_i > \al_j,
\end{cases}
\end{align}
and 
\begin{align}\label{def_B_ij}
B_{ij}(\alb,\beb) := \frac{1}{d_{ij}(\beb)}\begin{cases}
-\widetilde {P}_{i}^-(\alb,\beb), & \quad\text{ if }\ \al_i < \al_j,\\
\widetilde {P}_{i}^+(\alb,\beb), & \quad\text{ if }\ \al_i > \al_j,
\end{cases}
\end{align}
where 
\begin{align*}
\widetilde {P}_{i}^+(\alb,\beb)
 = \sum_{\substack{k=1, k\neq j}}^{\N}\max\{0, f_{ik}\},\quad\text{and}\quad\widetilde {P}_{i}^-(\alb,\beb) 
 = \sum_{\substack{k=1, k\neq j}}^{\N}\min\{0, f_{ik}\}.
\end{align*}
If  $d_{ij} = 0$, then we define $A_{ij}=B_{ij}=0$. 
Following {\normalfont{\cite[Lemma 3.5]{gabriel2016}}},  $A_{ij}$ and $B_{ij}$ are continuous functions and are Lipschitz 
with respect to the first variable $\alb$, i.e., there exist constants $\Lambda_{ij}^A(\beb,q_i)$ and 
$\Lambda_{ij}^B(\beb,q_i)$ such that

 \begin{align}
 |A_{ij}(\alb,\beb) - A_{ij}(\widetilde\alb,\beb)|&\le \Lambda_{ij}^A(\beb,q_i)\sum_{\ell\in\Nh(\omega_i)}|\al_\ell-\widetilde \al_\ell|,\label{Lipschitz_A_ij-0}\\
 |B_{ij}(\alb,\beb) - B_{ij}(\widetilde\alb,\beb)|&\le \Lambda_{ij}^B(\beb,q_i)\sum_{\ell\in\Nh(\omega_i)}|\al_\ell-\widetilde \al_\ell|.\label{Lipschitz_B_ij-0}
\end{align} 
 To show \eqref{Lipschitz_A_ij-0}--\eqref{Lipschitz_B_ij-0} we will employ  the following inequalities from \cite[Ineq. (2.3b)-(2.3c)]{lohmann2019}
\begin{align}\label{lohmann:ineq1}
|\max (a,b) - \max(\tilde{a}, \tilde{b})| &\leq |a - \tilde{a}| + |b - \tilde{b}|,\;\;\forall\,a,\tilde{a},b,\tilde{b}\in \mathbb{R},\\
|\min (a,b) - \min(\tilde{a}, \tilde{b})| &\leq |a - \tilde{a}| + |b - \tilde{b}|,\;\;\forall\,a,\tilde{a},b,\tilde{b}\in \mathbb{R}.
\label{lohmann:ineq2}
\end{align}

Let us assume $d_{ij}(\beb) > 0$ and
that  $\al_i - \al_j$ and $\widetilde\al_i - \widetilde\al_j$ are both positive. 
The other case for the differences $\al_i - \al_j<0$ and $\widetilde\al_i - \widetilde\al_j<0$ can be treated  analogously.
In view of \eqref{def_A_ij} and using \eqref{lohmann:ineq2}, we get
\begin{align}
|A_{ij}(\alb,\beb) - A_{ij}(\widetilde\alb,\beb)| & 
= \frac{1}{d_{ij}(\beb)}|\min\{ P_{i}^+(\alb,\beb), Q_{i}^+(\alb)\} - \min\{ P_{i}^+(\widetilde\alb,\beb), Q_{i}^+(\widetilde\alb)\}|\notag\\
& \le \frac{1}{d_{ij}(\beb)}( | P_{i}^+(\alb,\beb) -  P_{i}^+(\widetilde\alb,\beb)| + | Q_{i}^+(\alb) -  Q_{i}^+(\widetilde\alb)| ).\label{Lipschitz_A_ij}
\end{align}

Then, using  the definition of $P_i^+$, given by Algorithm \ref{algorithm-1}, we get

\begin{align}
&| P_{i}^+(\alb,\beb) -  P_{i}^+(\widetilde\alb,\beb)|  \leq \sum_{\ell\in \Nh}|\max\{0, f_{i\ell}(\alb,\beb)\} - \max\{0, f_{i\ell}(\widetilde\alb,\beb)\}|\notag\\
&\quad\leq \sum_{\ell\in \Nh,\,\ell\neq j}|\max\{0, f_{i\ell}(\alb,\beb)\} - \max\{0, f_{i\ell}(\widetilde\alb,\beb)\}|
 + |\max\{0, f_{ij}(\alb,\beb)\} - \max\{0, f_{ij}(\widetilde\alb,\beb)\}|\notag\\
&\quad
= I_1 + I_2,\label{P_i:bound}
\end{align}

In order to bound $I_1$ we will consider various cases for the the sign of the differences $\al_i - \al_\ell,\,\widetilde\al_i - \widetilde\al_\ell$, $\ell\in \Nh\;\ell\neq j.$ Let $\al_i - \al_\ell \leq 0$ and $\widetilde\al_i - \widetilde\al_\ell > 0, $ then
\begin{align}
&|\max\{0, f_{i\ell}(\alb,\beb)\} - \max\{0, f_{i\ell}(\widetilde\alb,\beb)\}| \leq d_{i\ell}(\beb)|\widetilde\al_i - \widetilde\al_\ell|  = d_{i\ell}(\beb)(\widetilde\al_i - \widetilde\al_\ell) 
 = d_{i\ell}(\beb)( (\widetilde\al_i - \al_i) + (\al_i - \widetilde\al_\ell) \notag\\
 &\qquad
 \leq d_{i\ell}(\beb)( (\widetilde\al_i - \al_i) + (\al_\ell - \widetilde\al_\ell)
 \le d_{i\ell}(\beb)( |\al_i - \widetilde\al_i| + |\al_\ell - \widetilde\al_\ell|).\label{l_1^l:1}
\end{align}
In a similar way we can show that for $\al_i - \al_\ell > 0$ and $\widetilde\al_i - \widetilde\al_\ell \leq 0,$ we get
\begin{align}\label{l_1^l:2}
|\max\{0, f_{i\ell}(\alb,\beb)\} - \max\{0, f_{i\ell}(\widetilde\alb,\beb)\}| \leq d_{i\ell}(\beb)|\al_i - \al_\ell| 
 \leq d_{i\ell}(\beb)( |\al_i - \widetilde\al_i| + |\al_\ell - \widetilde\al_\ell|).
\end{align}
Also if $\al_i - \al_\ell >0$ and $\widetilde\al_i - \widetilde\al_\ell > 0, $ then in view of \eqref{lohmann:ineq1}, we have
\begin{align}\label{l_1^l:3}
|\max\{0, f_{i\ell}(\alb,\beb)\} - \max\{0, f_{i\ell}(\widetilde\alb,\beb)\}|  \leq d_{i\ell}(\beb)( |\al_i - \widetilde\al_i| + |\al_\ell - \widetilde\al_\ell|),
\end{align}
Finally, for $\al_i - \al_\ell < 0$ and $\widetilde\al_i - \widetilde\al_\ell < 0$,
we get
\begin{align}\label{l_1^l:4}
|\max\{0, f_{i\ell}(\alb,\beb)\} - \max\{0, f_{i\ell}(\widetilde\alb,\beb)\}|  = 0.
\end{align}
Therefore, combining the above cases \eqref{l_1^l:1}--\eqref{l_1^l:4}, and the fact that $d_{i\ell}=0$ for $\ell\notin\Nh(\omega_i)$, we obtain 
\begin{equation}\label{eq:l_1}
I_1 \leq \sum_{\ell\in \Nh(\omega_i)}d_{i\ell}(\beb)( |\al_i - \widetilde\al_i| + |\al_\ell - \widetilde\al_\ell|).
\end{equation}

Moreover,
\begin{equation}\label{eq:l_2}
I_2 \leq d_{ij}(\beb)(|\al_i - \widetilde \al_i| + |\al_j - \widetilde\al_j|).
\end{equation}
Hence, using \eqref{eq:l_1} and \eqref{eq:l_2} in \eqref{P_i:bound}, we get

\begin{equation}\label{Lipschitz_p}
| P_{i}^+(\alb,\beb) -  P_{i}^+(\widetilde\alb,\beb)|  \leq C\max_{1\leq j\leq \N}d_{ij}(\beb)\sum_{j\in \Nh(\omega_i)}|\al_j - \widetilde\al_j| .
\end{equation}

Next, let $\al_i^{\max}$ and $\widetilde\al_i^{\max}$  the local maximum and local minimum, respectively, of  $\al_j$ and 
$\widetilde\al_j$, $j\in\Nh(\omega_i)$  for $\alb$ and $\widetilde\alb$, respectively, 
which may occur on different nodes. Therefore, let $\ell,\,\widetilde \ell\in \Nh(\omega_i)$ such that $\al_i^{\max} = \al_\ell$ and $\widetilde\al_i^{\max} = \widetilde\al_{\widetilde\ell}$, 
 then we obtain
 \begin{align}
| Q_{i}^+(\alb) -  Q_{i}^+(\widetilde\alb)| & = q_i|(\al_i^{\max} - \al_i) - (\widetilde\al_i^{\max} - \widetilde\al_i)| 
 \leq q_i|(\al_i^{\max} - \widetilde\al_i^{\max}) - (\al_i - \widetilde\al_i)|\notag\\
& \leq q_i|\al_i^{\max} - \widetilde\al_i^{\max}| + q_i|\al_i - \widetilde\al_i|= q_i|\al_\ell - \widetilde\al_{\widetilde\ell}| + q_i|\al_i - \widetilde\al_i|,\label{Q_i_bound}
\end{align}
We have 
\begin{equation}\label{al^max}
|\al_i^{\max} - \widetilde\al_i^{\max}|= \begin{cases}
\al_\ell - \widetilde\al_{\widetilde\ell}\le \al_\ell-\widetilde\al_\ell=| \al_\ell-\widetilde\al_\ell|,&\text{ if } \al_\ell>\widetilde\al_{\widetilde\ell},\text{ and }\ell\neq \widetilde\ell,\\
\widetilde\al_{\widetilde\ell} - \al_\ell \le  \widetilde\al_{\widetilde\ell} - \al_{\widetilde\ell}=|\al_{\widetilde\ell}-\widetilde\al_{\widetilde\ell}|,&\text{ if } \al_\ell<\widetilde\al_{\widetilde\ell},\text{ and } \ell\neq \widetilde\ell,\\
|\al_{\widetilde\ell} - \widetilde\al_{\widetilde\ell}|, &\text{ if }  \ell= \widetilde\ell.
\end{cases}
\end{equation}
Therefore, using \eqref{al^max} in \eqref{Q_i_bound} we get

\begin{equation}\label{Lipschitz_q}
| Q_{i}^+(\alb) -  Q_{i}^+(\widetilde\alb)|  \leq q_i\sum_{\ell\in \Nh(\omega_i)}|\al_\ell - \widetilde\al_\ell|.
\end{equation}
Hence, combining \eqref{Lipschitz_p} and \eqref{Lipschitz_q} in \eqref{Lipschitz_A_ij} we get 

\begin{align*}
|A_{ij}(\alb,\beb) - A_{ij}(\widetilde\alb,\beb)| \leq C\frac{\max_{1\leq j\leq \N}d_{ij}(\beb) + q_i}{d_{ij}(\beb)}\sum_{\ell\in \Nh(\omega_i)}|\al_\ell - \widetilde\al_\ell| .
\end{align*}

Similarly, in view of the definition of $B_{ij}$ we get, for $d_{ij}>0$,
\begin{align*}
|B_{ij}(\alb,\beb) - B_{ij}(\widetilde\alb,\beb)| & = \frac{1}{d_{ij}(\beb)}| \widetilde {P}_{i}^+(\alb,\beb) -  \widetilde {P}_{i}^+(\widetilde\alb,\beb)| \leq C\frac{\max_{1\leq j\leq \N}d_{ij}(\beb)}{d_{ij}(\beb)}\sum_{\ell\in \Nh(\omega_i)} |\al_\ell - \widetilde\al_\ell|.
\end{align*}
 \end{proof}

\begin{definition} Let $w\in\Sh$, then
$$\Lambda_{ij}(w,q_i) := (d_{ij}^{-1}(w)(\max_{1\leq j\leq \N}d_{ij}(w) + q_i) + 1).$$
\end{definition}

\begin{lemma}\label{lemma:limiters_estimate_localy}
Let  $\chi,\psi,w\in \Sh$,  
and $\widetilde {\rho}_{ij}(\chi,w):=\rho_{ij}(\chi,w)(\chi_i - \chi_j)$,  with $i,j\in\Nh$, and $\rho_{ij}$ given in \eqref{stab_term_afc_2D}.  Then there exists $C>0$ independent of $h$ such that $\widetilde {\rho}_{ij}$ satisfies,
\begin{align*}
\vert \widetilde {\rho}_{ij}(\chi,w) - \widetilde {\rho}_{ij}(\psi,w)\vert \leq C\Lambda_{ij}(w,q_i)\sum_{\ell\in \Nh(\omega_i)}|\chi_\ell - \psi_\ell|, \quad \forall \chi,\psi,w\in \Sh.
\end{align*}
\end{lemma}
\begin{proof} It suffices to consider the case where $\rho_{ij}=1 - \mathfrak{a}_{ij}$.
We easily get
\begin{align*}
\vert\widetilde {\rho}_{ij}(\chi,w) - \widetilde {\rho}_{ij}(\psi,w)\vert & =  \vert (1 - \mathfrak{a}_{ij}(\chi,w))(\chi_i - \chi_j) - (1 - \mathfrak{a}_{ij}(\psi,w))(\psi_i - \psi_j) \vert \\
& \leq \vert (\chi_i - \chi_j) - (\psi_i - \psi_j)\vert + \vert\Phi_{ij}(\chi,w) - \Phi_{ij}(\psi,w)\vert,
\end{align*}
where  
$\Phi_{ij}(\chi,w) := \mathfrak{a}_{ij}(\chi,w)(\chi_i - \chi_j)$, $\chi\in\Sh$, $j\in \Nh$.
 
Let $\chi,\psi\in \Sh,$ then following of the proof of \cite[Lemma 3.5]{gabriel2016}, we have for 
$(\chi_i - \chi_j)(\psi_i - \psi_j) \leq 0$ that
\begin{align*}
|\Phi_{ij}(\chi,w) - \Phi_{ij}(\psi,w)| & \leq |\chi_i - \chi_j| + |\psi_i - \psi_j| 
 = |(\chi_i - \psi_i) - (\chi_j - \psi_j)| 
 \leq C\sum_{l\in \Nh(\omega_i)}|\chi_l - \psi_l|.
\end{align*}

Further, again in view of the proof of \cite[Lemma 3.5]{gabriel2016}, we have for $(\chi_i - \chi_j)(\psi_i - \psi_j) > 0$ that
\begin{align*}
\Phi_{ij}(\chi,w) - \Phi_{ij}(\psi,w) &=(A_{ij}(\bchi,\bw)-A_{ij}(\bpsi,\bw))\dfrac{\psi_i - \psi_j}{|\psi_i - \psi_j|+B_{ij}(\bpsi,\bw)}\\
&+\mathfrak{a}_{ij}(\bchi,\bw)\dfrac{(B_{ij}(\bpsi,\bw)-B_{ij}(\bchi,\bw))(\psi_i - \psi_j)+((\chi_i - \chi_j)-(\psi_i - \psi_j))B_{ij}(\bpsi,\bw)}{|\psi_i - \psi_j|+B_{ij}(\bpsi,\bw)},
\end{align*}
where $\bchi=(\chi_1,\dots,\chi_\N)^T$ and $\bw=(w_1,\dots,w_\N)^T$.
 Therefore the desired bound follows from
\begin{align*}
|\Phi_{ij}(\bchi,\bw) - \Phi_{ij}(\bpsi,\bw)| & \leq |A_{ij}(\bchi,\bw) - A_{ij}(\bpsi,\bw)| 
+ |B_{ij}(\bchi,\bw) - B_{ij}(\bpsi,\bw)| + |(\chi_i - \psi_i) - (\chi_j - \psi_j)|\\
&\le (\Lambda_{ij}^A+\Lambda_{ij}^B)\sum_{\ell\in\Nh(\omega_i)}|\chi_\ell - \psi_\ell|+ |(\chi_i - \psi_i) - (\chi_j - \psi_j)|.
\end{align*}
\end{proof}

\begin{remark}\label{remark:L_bound} 
The Lipschitz constants $\Lambda_{ij}$, in Lemma \ref{lemma:limiters_estimate_localy}  can be estimated in view of \eqref{est_dij_2D_L_infty}, as follows,
\begin{align*}
d_{ij}(w)\Lambda_{ij}(w,q_i) &= \max_{1\leq j\leq \N}d_{ij}(w) + q_i
\leq C(h\,\|\nabla w\|_{L^{\infty}} + \|\q\|_{\max}), 
\end{align*}
where $\|\q\|_{\max} = \max_{i\in\Nh}|q_i|$, $\q=(q_1,\dots,q_N)$ and the coefficients $q_i$, $i\in \Nh$,  in Algorithm \ref{algorithm-1}, are such that Lemma \ref{corrollary:linearity_preservation_limiters} holds. 
\end{remark}

\begin{lemma}\label{lemma:limiters_estimate_localy-ver2}
Let  $\chi,\psi,w\in \Sh$,  
and $\widehat {\rho}_{ij}(\chi,w):=d_{ij}(w)\rho_{ij}(\chi,w)(\chi_i - \chi_j)$,  with $i,j\in\Nh$ and $\rho_{ij}$ given in \eqref{stab_term_afc_2D}.  Then there exists $C>0$ independent of $h$ such that $\widehat {\rho}_{ij}$ satisfies,
\begin{align*}
\vert \widehat {\rho}_{ij}(\chi,w) - \widehat {\rho}_{ij}(\chi,\widehat w)\vert \leq 
Ch\|\nabla (w-\widehat w)\|_{L^{\infty}(\omega_i)} \|\nabla \chi\|_{L^{2}(\omega_i)}
\end{align*}
\end{lemma}
\begin{proof} It suffices to consider the case where $\rho_{ij}=1 - \mathfrak{a}_{ij}$.
We easily get
\begin{align}\label{eq:I-lemma_last-rho}
\vert\widehat {\rho}_{ij}(\chi,w) - \widehat {\rho}_{ij}(\chi,\widehat w)\vert & =  \vert d_{ij}(w)(1 - \mathfrak{a}_{ij}(\chi,w))(\chi_i - \chi_j) - d_{ij}(\widehat w)(1 - \mathfrak{a}_{ij}(\chi,\widehat w))(\chi_i - \chi_j) \vert \notag\\
& \leq \vert (d_{ij}(w)-d_{ij}(\widehat w))(\chi_i - \chi_j) \vert + \vert\Phi_{ij}(\chi,w) - \Phi_{ij}(\chi,\widehat w)\vert,
\end{align}
where $\Phi_{ij}(\chi,w) := d_{ij}(w)\mathfrak{a}_{ij}(\chi,w)(\chi_i - \chi_j)$, $\chi\in\Sh$, $j\in \Nh$.

Following the proof in Lemma \ref{lemma:limiters_estimate_localy} we have for $\chi_i - \chi_j> 0$ that
\begin{align}\label{eq:I-lemma_last-Phi}
|\Phi_{ij}(\chi,w) - \Phi_{ij}(\chi,\widehat w)| &=|d_{ij}(w)A_{ij}(\chi,w)\dfrac{\chi_i - \chi_j}{|\chi_i - \chi_j|+B_{ij}(\chi, w)} 
-d_{ij}(\widehat w)A_{ij}(\chi,\widehat w)\dfrac{\chi_i - \chi_j}{|\chi_i - \chi_j|+B_{ij}(\chi,\widehat w)}|\notag\\
&\le|d_{ij}(w)A_{ij}(\chi,w)-d_{ij}(\widehat w)A_{ij}(\chi,\widehat w)||\dfrac{\chi_i - \chi_j}{|\chi_i - \chi_j|+B_{ij}(\chi, w)}|\notag\\
&\quad+|d_{ij}(\widehat w)A_{ij}(\chi,\widehat w)(\dfrac{\chi_i - \chi_j}{|\chi_i - \chi_j|
+B_{ij}(\chi,w)}-\dfrac{\chi_i - \chi_j}{|\chi_i - \chi_j|+B_{ij}(\chi,\widehat w)})|
=I+II.
\end{align}

Let us assume $d_{ij}(w) > 0$, $d_{ij}(\widehat w) > 0$, and  $\chi_i - \chi_j>0$.
The other case for the differences $\chi_i - \chi_j<0$  can be treated  analogously.
In view of \eqref{def_A_ij} and using \eqref{lohmann:ineq2}, we get
\begin{align}
I\le |d_{ij}(w)A_{ij}(\chi,w)-d_{ij}(\widehat w)A_{ij}(\chi,\widehat w)| & 
= |\min\{ P_{i}^+(\chi,w), Q_{i}^+(\chi)\} - \min\{ P_{i}^+(\chi,\widehat w), Q_{i}^+(\chi)\}|\notag\\
& \le  | P_{i}^+(\chi,w) -  P_{i}^+(\chi,\widehat w)|.\label{Lipschitz_A_ij-1}
\end{align}
 
Then, using  the definition of $P_i^+$, given by Algorithm \ref{algorithm-1}, we get

\begin{align}
&|  P_{i}^+(\chi,w) -  P_{i}^+(\chi,\widehat w)|  \leq \sum_{\ell}|\max\{0, f_{i\ell}(\chi,w)\} - \max\{0, f_{i\ell}(\chi,\widehat w)\}|.
\label{P_i:bound-1}
\end{align}

In order to the right hand-side of \eqref{P_i:bound-1} we will consider various cases for the the sign of the difference $\chi_i - \chi_\ell$, 
$\ell\in \Nh\;\ell\neq i.$ Let $\chi_i - \chi_\ell > 0,$ then
\begin{align}
|\max\{0, f_{i\ell}(\chi,w)\} &- \max\{0, f_{i\ell}(\chi,\widehat w)\}| \leq |d_{i\ell}(w)-d_{i\ell}(\widehat w)||\chi_i - \chi_\ell|.\label{l_1^l:1-1}
\end{align}
For $\chi_i - \chi_\ell \le 0$
we get
\begin{align}\label{l_1^l:2-1}
|\max\{0, f_{i\ell}(\chi,w)\} - \max\{0, f_{i\ell}(\chi,\widehat w)\}|  = 0.
\end{align}
Therefore, combining the above cases \eqref{l_1^l:1-1}--\eqref{l_1^l:2-1}, \eqref{P_i:bound-1}, \eqref{diff-tau_ij} and the fact that $d_{i\ell}=0$ for $\ell\notin\Nh(\omega_i)$, we obtain 
\begin{equation}\label{eq:l_1-1}
|  P_{i}^+(\chi,w) -  P_{i}^+(\chi,\widehat w)| \leq \sum_{l\in \Nh(\omega_i)}|d_{i\ell}(w)-d_{i\ell}(\widehat w)||\chi_i - \chi_j|
\le Ch\|\nabla (w-\widehat w)\|_{L^{\infty}(\omega_i)} \|\nabla \chi\|_{L^{2}(\omega_i)}.
\end{equation}

Hence, combining \eqref{Lipschitz_A_ij-1} and  \eqref{eq:l_1-1}  we get 
\begin{align}\label{lips_Aij}
I\le |d_{ij}(w)A_{ij}(\chi,w)-d_{ij}(\widehat w)A_{ij}(\chi,\widehat w)| \leq Ch\|\nabla (w-\widehat w)\|_{L^{\infty}(\omega_i)} \|\nabla \chi\|_{L^{2}(\omega_i)} .
\end{align}

We turn now to the estimation of $II$. In view of \eqref{correction_factors_formulation} we can rewrite $II$ as,
\begin{align}\label{eq:II-est}
&II= |d_{ij}(\widehat w)\dfrac{A_{ij}(\chi,\widehat w)}{|\chi_i - \chi_j|+B_{ij}(\chi,\widehat w)} \dfrac{\chi_i - \chi_j}{|\chi_i - \chi_j|+B_{ij}(\chi,w)}(B_{ij}(\chi,\widehat w)-B_{ij}(\chi, w))|\notag\\
&\le |\overline{\mathfrak{a}}_{ij}(x,\widehat w)\dfrac{\chi_i - \chi_j}{|\chi_i - \chi_j|+B_{ij}(\chi,w)}
\{(d_{ij}(\widehat w)B_{ij}(\chi, \widehat w)-d_{ij}(w)B_{ij}(\chi, w))+(d_{ij}( w)-d_{ij}(\widehat w))B_{ij}(\chi, w)\}|\notag\\
&=|\overline{\mathfrak{a}}_{ij}(x,\widehat w)\dfrac{\chi_i - \chi_j}{|\chi_i - \chi_j|+B_{ij}(\chi,w)}|
|d_{ij}(\widehat w)B_{ij}(\chi, \widehat w)-d_{ij}(w)B_{ij}(\chi, w)|\notag\\
&+|\overline{\mathfrak{a}}_{ij}(x,\widehat w)\dfrac{\chi_i - \chi_j}{|\chi_i - \chi_j|+B_{ij}(\chi,w)}||(d_{ij}( w)-d_{ij}(\widehat w))B_{ij}(\chi, w)|\notag\\
&\le |d_{ij}(\widehat w)B_{ij}(\chi, \widehat w)-d_{ij}(w)B_{ij}(\chi, w)|
+|\chi_i - \chi_j||d_{ij}( w)-d_{ij}(\widehat w)|\le II_1+II_2
\end{align}
Similarly as we derived \eqref{lips_Aij} and in view of the definition of $B_{ij}$ we get, for $d_{ij}>0$,
\begin{align}\label{eq:II-est-2}
II_1\le   | \widetilde {P}_{i}^+(\chi,w) -  \widetilde {P}_{i}^+(\chi,\widehat w)|\le Ch\|\nabla (w-\widehat w)\|_{L^{\infty}(\omega_i)} \|\nabla \chi\|_{L^{2}(\omega_i)}.
\end{align}
Further
\begin{align}\label{eq:II-est-3}
II_2=
 |d_{ij}( w)-d_{ij}(\widehat w)| |\chi_i - \chi_j|
&\le Ch\|\nabla (w-\widehat w)\|_{L^{\infty}(\omega_i)} \|\nabla \chi\|_{L^{2}(\omega_i)}.
\end{align}
Therefore \eqref{eq:II-est}, \eqref{eq:II-est-2} and \eqref{eq:II-est-3} give
\begin{align}\label{eq:II-est-final}
 |II|\le C\|\nabla (w-\widehat w)\|_{L^{2}(\omega_i)} \|\nabla \chi\|_{L^{2}(\omega_i)}
\end{align}
Hence combining \eqref{eq:I-lemma_last-rho}, \eqref{eq:I-lemma_last-Phi}, \eqref{lips_Aij}, 
\eqref{eq:II-est-final} we obtain the derised result.
\end{proof}

In the sequel we will demonstrate Lemmas \ref{lemma:estimate_stab_term}-\ref{lemma:chemotaxis_bounds_2D-1}.
\begin{proof}[Proof of Lemma \ref{lemma:estimate_stab_term}]
 
Following  the proof of Lemma \ref{lemma:estimate_1_d_h_afc} and  \eqref{diff-tau_ij}, we obtain the desired bound
\begin{align*}
|d_h(w;\psi, \chi) - d_h(\widetilde {w};\psi,\chi)| 
& =| \sum_{i<j}(d_{ij}(w) - d_{ij}(\widetilde {w}))(\psi_i - \psi_j)(\chi_i-\chi_j)|\\
&\leq Ch\|\nabla (w - \widetilde {w})\|_{L^{\infty}}\sum_{i<j}\vert \psi_i - \psi_j\vert\vert\chi_i - \chi_j\vert
 \leq Ch\|\nabla (w - \widetilde {w})\|_{L^{\infty}}\|\nabla \psi\|\|\nabla \chi\|.
\end{align*}
\end{proof}

\begin{proof}[Proof of Lemma \ref{lemma:estimate_1_d_h_afc-ver2}]
Using \eqref{equiv_stab_term_afc}, we have, 
\begin{equation}
\overline{d}_h(\psi,w;\psi, \chi) - \overline{d}_h(\psi,\widehat w;\psi, \chi)  =
 \sum_{i<j}(\widehat {\rho}_{ij}(\psi,w) - \widehat {\rho}_{ij}(\psi,\widehat w))(\chi_i - \chi_j),\nonumber
\end{equation}
where $\widehat {\rho}_{ij}(\psi,w) : = d_{ij}(w)\rho_{ij}(\psi,w)(\psi_i - \psi_j),\;\forall \psi,w\in \Sh$. 
In view of Lemma \ref{lemma:limiters_estimate_localy-ver2}, 
we get the desired result
\begin{align*}
\vert\overline{d}_h(\psi,w;\psi, \chi) - \overline{d}_h(\psi,\widehat w;\psi, \chi)\vert & \leq C \sum_{i<j}\vert\widehat {\rho}_{ij}(\psi, w) - \widehat {\rho}_{ij}(\psi,\widehat w)\vert \vert\chi_i - \chi_j\vert
 \leq Ch\|\nabla (w-\widehat w)\|_{L^{\infty}} \|\nabla \psi\|\|\nabla \chi\|.
\end{align*}
\end{proof}

\begin{proof}[Proof of Lemma \ref{lemma:stability_stab_term2}]
We will follow the proof of \cite[Lemma 3]{gabriel2018}.
Using \eqref{equiv_stab_term_afc}, we have
\begin{align}
\overline{d}_h(\psi,w;\psi,\chi) & = \sum_{i<j}d_{ij}(w){\rho}_{{ij}}(\psi,w)(\psi_i - \psi_j)(\chi_i - \chi_j)=\sum_{i<j}d_{ij}(w)\widetilde {\rho}_{ij}(\psi,w)(\chi_i - \chi_j)
=I,\label{lemma:stability_eq1}
\end{align}
where $\widetilde {\rho}_{ij}(\psi,w)={\rho}_{{ij}}(\psi,w)(\psi_i - \psi_j)$.

Given $i\in\Nh$, let  $\widehat{v}_i\in {\P}_1(\omega_{i})$ be the unique solution of 
\begin{align*}
(\nabla \widehat{v}_{i}, \nabla \chi)_{L^{2}(\omega_{i})} & = (\nabla v, \nabla \chi)_{L^{2}(\omega_{i})},\quad\forall\,\chi\in {\P}_1(\omega_{i}),\\
(\widehat{v}_{i},1)_{L^{2}(\omega_{i})} & = (v,1)_{L^{2}(\omega_{i})}.
\end{align*}
Also, in view of  Assumption \ref{mesh-assumption} there exists a constant $C$ independent of $h$, such that
\begin{equation}\label{eq:local-estimate}
\|\nabla (v - \widehat{v}_{i})\|_{L^{2}(\omega_{i})} \leq Ch\|v\|_{H^{2}(\omega_{i})}.
\end{equation}
We can extended   $\widehat{v}_{i}$  as a function  on  ${\P}_1(\omega_{i})$. In view of Lemma 
\ref{corrollary:linearity_preservation_limiters} and the definition of the linearity preservation, i.e., 
\eqref{eqn:linear_preserve}, we have that 
$\rho_{ij}(\widehat{v}_{i},w) = 0$, for $j\in\Nh$. 

Thus $I$ in \eqref{lemma:stability_eq1} can be rewritten as
\begin{equation}\label{d_h_est_2D-0}
\begin{aligned}
I&=\sum_{i<j }d_{ij}(w)\widetilde{\rho}_{ij}(\psi,w)(\chi_i - \chi_j)
  =  \sum_{i<j} d_{ij}(w) (\widetilde{\rho}_{ij}(\psi,w)  - \widetilde{\rho}_{ij}(\widehat{v}_{i},w) )(\chi_i - \chi_j).
\end{aligned}
\end{equation}

Then, we easily get
\begin{align*}
 |I|&  \leq \sum_{i<j}d_{ij}(w)| \widetilde{\rho}_{ij}(\psi,w) - \widetilde{\rho}_{ij}(\widehat{v}_{i},w)| |\chi_i - \chi_j|
 \leq \Bigl(\sum_{i<j} d_{ij}(w)^2|\widetilde{\rho}_{ij}(\psi,w) - \widetilde{\rho}_{ij}(\widehat{v}_{i},w)|^2\Bigr)^{1/2}\bigl(\sum_{i ,j= 1}^{\N}|\chi_i - \chi_j|^2\bigr)^{1/2}\\
& \qquad \leq C \Bigl(\sum_{i<j}d_{ij}(w)^2|\widetilde{\rho}_{ij}(\psi,w) - \widetilde{\rho}_{ij}(\widehat{v}_{i},w)|^2\Bigr)^{1/2}\|\nabla \chi\|.
\end{align*}

In view of  Lemma \ref{lemma:limiters_estimate_localy}, \eqref{eq:local-estimate} and  Remark \ref{remark:L_bound} we obtain

\begin{align*}
 d_{ij}(w) &|\widetilde{\rho}_{ij}(\psi,w)  - \widetilde{\rho}_{ij}(\widehat{v}_{i},w)|
  \leq C d_{ij}(w) \Lambda_{ij}(w,q_i) (\sum_{l\in \Nh(\omega_i)} |\psi_l - (\widehat{v}_{i})_l|)^{1/2}\\
& \leq C(h\|\nabla w\|_{L^{\infty}} + \|\q\|_{\max}) \|\nabla (\psi - \widehat{v}_{i})\|_{L^{2}(\omega_i)} \\
& \leq C(h\|\nabla w\|_{L^{\infty}} + \|\q\|_{\max}) (\|\nabla (\psi - v)\|_{L^{2}(\omega_{i})} 
+ \|\nabla (v - \widehat{v}_{i})\|_{L^{2}(\omega_{i})})\\
 &\leq C(h\|\nabla w\|_{L^{\infty}} + \|\q\|_{\max}) (\|\nabla (\psi - v)\|_{L^{2}(\omega_{i})} + h\|v\|_{H^{2}(\omega_{i})}),
\end{align*}
which gives the desired result
\begin{equation*}
|I|\le C(h\,\|\nabla w\|_{L^{\infty}} + \|\q\|_{\max})(\|\nabla (\psi - v)\| + h\|v\|_{2} )\|\nabla\chi\|.
\end{equation*}
\end{proof}

\begin{proof}[Proof of Lemma \ref{lemma:estimate_2_d_h_afc}]
Using \eqref{equiv_stab_term_afc}, we have, 
\begin{equation}
\overline{d}_h(v,w; v, \chi) - \overline{d}_h(s,w; s, \chi)  =
 \sum_{i<j}d_{ij}(w)(\widetilde {\rho}_{ij}(v,w) - \widetilde {\rho}_{ij}(s,w))(\chi_i - \chi_j).\nonumber
\end{equation}
where $\widetilde {\rho}_{ij}$ are defined in Lemma \ref{lemma:limiters_estimate_localy}.
In view of Remark \ref{remark:L_bound} and Lemma \ref{lemma:limiters_estimate_localy}, we get the desired estimate.
\end{proof}

\begin{proof}[Proof of Lemma \ref{lemma:chemotaxis_bounds_2D-1}]
We easily get the following splitting
\begin{align}
(u\nabla c - v\nabla w, \nabla \chi) & =  (u\nabla (c - R_hc), \nabla \chi)
 + (u\nabla (R_hc -  w), \nabla \chi)
 + ((u - R_hu)\nabla w, \nabla \chi)
 + ((R_hu -  v)\nabla w, \nabla \chi)\notag\\
& = I_1 + \dots + I_4.\nonumber
\end{align}
Next, we will bound each one of the terms $I_i$, $i=1,2,3,4$. Note that  by applying integration by parts $I_1$ can be rewritten  as
\begin{equation}\label{I_1}
I_1 =    (u\nabla \chi, \nabla (c - R_hc))= \sum_{K\in\Th} \left( \int_K\nabla \cdot (u \nabla \chi)(c - R_hc) \,dx - \int_{\partial K} (c - R_hc) u\nabla \chi  \cdot n_K\,ds \right),
\end{equation}
where $n_K$ is the unit normal vector on $\partial K.$ 
 Employing now the  trace inequality \eqref{eq:scaled_trace}, the error estimates for $R_h$, \eqref{ritz_projection_inf}, and the fact that $\chi$ is linear on $K$, we obtain
\begin{align*}
\left\vert\int_{\partial K}  (c - R_hc) u\nabla \chi  \cdot n_K\,ds \right\vert & \le C\|c - R_hc\|_{L^{4}(\partial K)} \| u\|_{L^{4}(\partial K)} \|\nabla \chi\|_{L^{2}(\partial K)}  \\
& \le C\,h_K^{-1/2}h_K^{1/4}\|c - R_hc\|_{L^{\infty}(K)}\| u\|_{L^{4}(\partial K)}\|\nabla \chi\|_{L^{2}(K)} \\
& \le C\,h_K^{-1/2}h_K^{1/4}\|c - R_hc\|_{L^{\infty}}\| u\|_{L^{4}(\partial K)}\|\nabla \chi\|_{L^{2}(K)} \\
& \leq C\,h^{3/2+1/4}|\log h|\|c\|_{2,\infty}\| u\|_{L^{4}(\partial K)}\|\nabla \chi\|_{L^{2}(K)}.
\end{align*}
We note that in view of the Gagliardo--Nirenberg--Ladyzhenskaya inequality, we get
\begin{align*}
\|u\|_{L^{4}} \leq C_{GNL}(\Omega)\|u\|_{L^{2}}^{1/2}\|\nabla u\|_{L^{2}}^{1/2}.
\end{align*}
Then, we perform a scaling argument, using e.g., \cite[Eq. (4.5.2)]{brenner2008}, to obtain a constant $C,$ independent of $K\in\Th,$ such that 
\begin{align*}
\|u\|_{L^{4}(K)} \leq C\|u\|_{H^{1}(K)},\;\;K\in\Th.
\end{align*}
To estimate the $\|u\|_{L^{4}(\partial K)},$ we use the trace inequality \eqref{eq:scaled_trace} for $p=4,$ and the latter estimate, to get
\begin{align*}
\|u\|_{L^{4}(\partial K)} \leq C\left(h_K^{-1/4}\|u\|_{H^{1}(K)} + h_K^{3/4}\|\nabla u\|_{L^{4}(K)}\right).
\end{align*}
Then, using the above estimates, yields
\begin{align*}
\left\vert\int_{\partial K}  (c - R_hc) u\nabla \chi  \cdot n_K\,ds \right\vert &  \leq  C(c)\,h^{3/2+1/4}|\log h|\left(h_K^{-1/4}\|u\|_{H^{1}(K)} + h_K^{3/4}\|\nabla u\|_{L^{4}(K)}\right)\|\nabla \chi\|_{L^{2}(K)}.
\end{align*}
Summing over the elements $K\in\Th,$ and using the fact that $\Th$ is quasi-uniform, we obtain
\begin{align*}
\sum_{K\in\Th}\left\vert  \int_{\partial K}   (c - R_hc) (u\nabla \chi) \cdot n\,ds \right\vert & \leq  C(c)\,h^{3/2}|\log h| \left(\sum_{K\in\Th} \|u\|_{H^{1}(K)}^2\right)^{1/2}\left(\sum_{K\in\Th}\|\nabla \chi\|^2_{L^{2}(K)}\right)^{1/2}\\
& +  C(u,c)\,h^{5/2}|\log h|\left(\sum_{K\in\Th} 1^2\right)^{1/2}\left(\sum_{K\in\Th}\|\nabla \chi\|_{L^{2}(K)}^2\right)^{1/2}\\
& \leq C(u,c)\,h^{3/2}|\log h|\|\nabla \chi\|.
\end{align*}
Further, since $\chi$ is linear on every $K\in\Th,$
\begin{align*}
\left\vert \int_K\nabla \cdot (u \nabla \chi)(c - R_hc) \,dx\right\vert & = \left\vert \int_K \nabla u \cdot \nabla \chi (c - R_hc) \,dx\right\vert \leq \|\nabla u\|_{L^{\infty}(K)}\|\nabla \chi\|_{L^{2}(K)}\|c - R_hc\|_{L^{2}(K)}.
\end{align*}
Summing over the elements $K\in\Th,$ 
\begin{align*}
\sum_{K\in\Th} \left\vert \int_K\nabla \cdot (u \nabla \chi)(c - R_hc) \,dx\right\vert & \leq C(u)\left(\sum_{K\in\Th}\|c - R_hc\|_{L^{2}(K)}^2\right)^{1/2}\left(\sum_{K\in\Th}\|\nabla \chi\|_{L^{2}(K)}^2\right)^{1/2}\\
& \leq C(u)\|c - R_hc\| \|\nabla \chi\| \leq C(u,c)h^2\|\nabla \chi\|,
\end{align*}
where we have used \eqref{ritz_projection_est2_2D}. Thus, we get for the term $I_1,$ that
\begin{equation*}
|I_1| \leq C(u,c)h^{3/2}|\log h|\|\nabla\chi\|,
\end{equation*} 
with constant $C(u,c)$ depends on norms of $u$ and $c,$ but independent of $h.$
Next, using again \eqref{ritz_projection_est2_2D}, we can easily bound $I_2, I_3$ and $I_4$ in the following way,
\begin{align*}
 |I_3|  &\leq C\|u - R_hu\|\|\nabla w\|_{L^{\infty}}\|\nabla \chi\| 
 \leq Ch^2\|\nabla w\|_{L^{\infty}}\|u\|_{2}\|\nabla \chi\|,\\
|I_2| &\leq \|u\|_{L^{\infty}}\|\nabla (w - R_hc)\| \|\nabla \chi\|,
\qquad |I_4|  \leq \|\nabla w\|_{L^{\infty}}\|v - R_hu\|\|\nabla \chi\|.
\end{align*}
Therefore, combining the above bounds for $I_i$, $i=1,2,3,4$, we obtain the desired estimate.
\end{proof}

\section{Numerical experiments}\label{section:numerical_results}
\begin{figure}
\centering
\includegraphics[scale=.8]{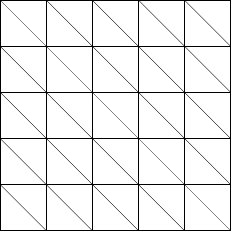}
\caption{A triangulation of a square domain.}\label{fig:triangulation}
\end{figure}

In this section we present several numerical experiments, illustrating our theoretical results. 
We consider a uniform mesh $\Th$ of the unit square $\Omega = [0,1]^2.$ Each side of $\Omega$ is divided into $M_0$ 
intervals of length $h_0=1/M_0$ for $M_0\in\mathbb{N}$ and we define the triangulation $\Th$ by dividing each small square 
by its diagonal, see Fig. \ref{fig:triangulation}. Thus $\Th$ consists of $2M_0^2$ right-angle triangles with diameter $h = \sqrt{2}h_0.$ Obviously $\Th$ satisfies Assumption \ref{mesh-assumption}. Therefore, the corresponding stiffness matrix $\S$ has non-positive off-diagonal elements and positive diagonal elements.  
To construct the approximation at each time level $t^n$ we implement the fixed point iteration scheme 
\eqref{dasdasdasda}. As a stopping criterion we consider the relative 
error between two successive solutions of \eqref{dasdasdasda}, in the maximum norm of $\mathbb{R}^{\N}$,
 with $\mathsf{TOL}=10^{-5}$. 
For the computation of the correction factors $\mathfrak{a}_{ij},$ we use Algorithm \ref{algorithm-1}
 with  $q_{i} = \sum_{j\neq i}d_{ij}$, where $d_{ij}$, $i,j=1,\dots,\N$, are the elements of $\D$. Due to the construction of $\Th$, $\gamma_i = 1,\; i=1,\ldots,\N$ in Remark \ref{remark:linearity_preservation}, see \cite{gabriel2017b}.

\subsection{Example with known solution}

\begin{figure}
\begin{tabular}{cc}
\includegraphics[scale=0.3]{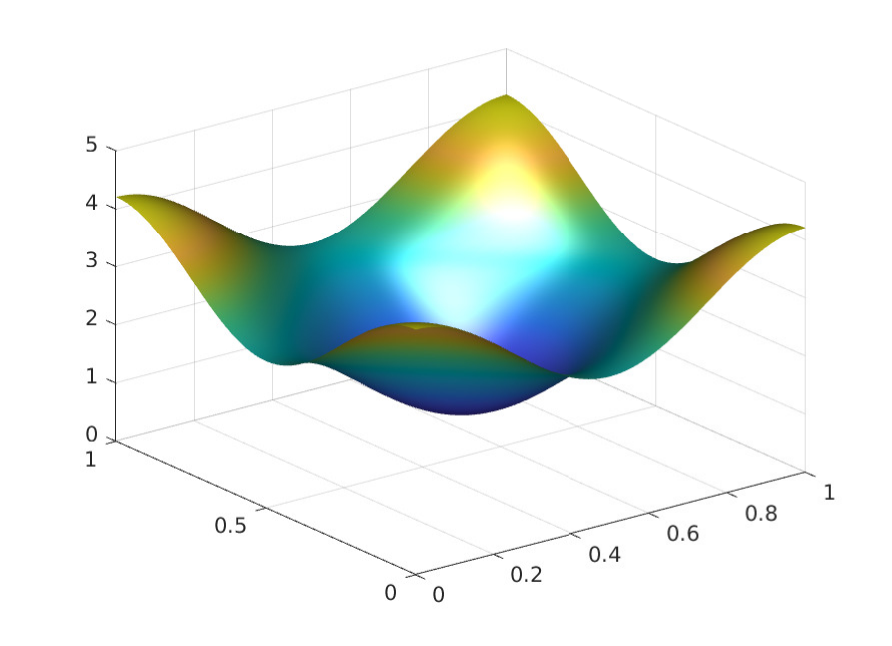} 
\includegraphics[scale=0.3]{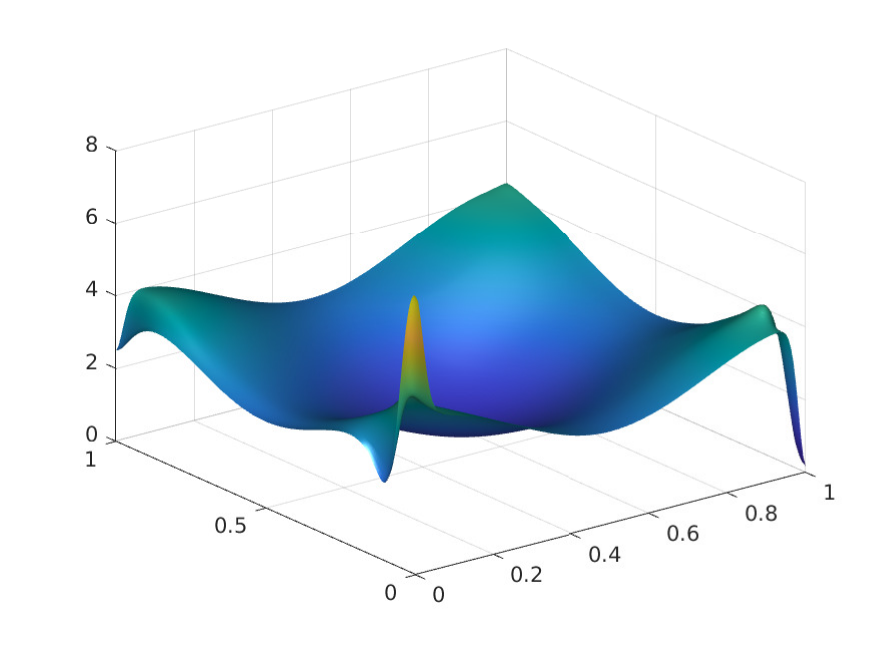}  
\includegraphics[scale=0.3]{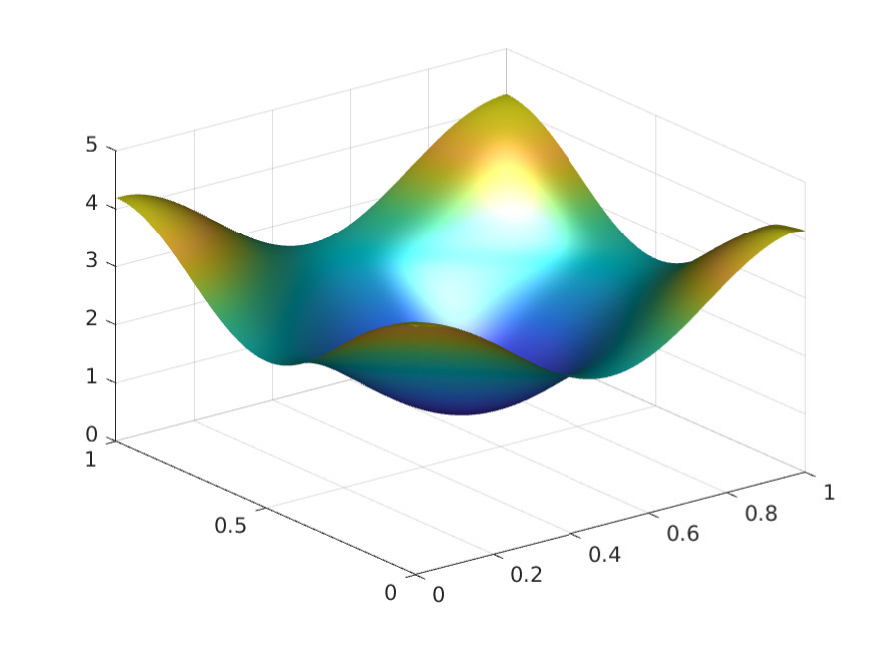}  
\end{tabular}
\caption{From left to right: Profile of the exact solution $u$ of \eqref{steep_gradients} at $T=1$, the low-order approximation, the AFC approximation.}\label{fig:profile_of_solution}
\centering
\end{figure}

In this subsection, we consider for $\lambda = \eta = 1,$  the source terms $f_u,\,f_c,$ on the \eqref{Minimal_model_uc}, so that the solution is given  for $(\x ,t)\in [0,1]^2 \times [0,1],$
\begin{equation}\label{steep_gradients}
\begin{aligned}
u(x,y,t) & = \sin(t)\left(\cos(2\pi x) + \cos(2\pi y) + 3\right),\\
c(x,y,t) & = \sin(t)(\cos(2\pi x) + \sin(2\pi y) - 2\pi y + 9).
\end{aligned}
\end{equation}
Note that $\frac{\partial u}{\partial \n} = \frac{\partial c}{\partial \n}  = 0$ on $\partial\Omega.$ In Fig. \ref{fig:profile_of_solution} with plot the exact solution $u$ \eqref{steep_gradients} at $T = 1$ with $\mu=10^{-2},$ together with the solution of the low-order and AFC scheme, respectively. For the approximations, we used $h_0 = 1/160$, corresponding to $25921$ degrees of freedom, and $k = 1/8100$. The low-order scheme fails to approximate the solution accurately; in particular, it produces large values at the four corners of the square. On the other hand, the AFC scheme appears to approximate the solution well, as will be confirmed in the forthcoming convergence analysis.

For the above example, we set $\mu = 1, 10^{-2},$ and we compute the experimental order of convergence of the low-order scheme as well as the AFC scheme in $L^{2}$ and $H^{1}$ norm. In particular, we compute the errors on a sequence of triangulations with $h_0 = 1/M_0,$ with $M_0=10,20,\ldots,160.$ The time step is $k = c_1h_0^{\log_{2} 3}\approx c_1h_0^{3/2},$ where $c_1 \approx \frac{100}{0.026}.$ Since also $k = T/N_0,$ this means that $N_0 = 100, 300, \ldots, 8100.$


\begin{table}
\begin{center}
\caption{Error and experimental convergence order for \eqref{steep_gradients} with $\mu = \eta = 1$.}
\label{table:man_convergenceL2_m00}
\begin{tabular}{@{}ccccc|cccc@{}}
\toprule
 & \multicolumn{4}{c}{$\max_{1\leq n\leq \NT}\|U^n-u^n\|$} & \multicolumn{4}{c}{$\max_{1\leq n\leq \NT}\|C^n-c^n\|_1$} \\ \midrule
$h_0$ & LOW & Order & AFC & Order & LOW & Order & AFC & Order \\ \midrule
$1/10$    & $2.0100e-1$ &       & $4.8261e-2$ &       & $9.4753e-1$ &       & $9.4718e-1$ &       \\ 
$1/20$    & $1.0482e-1$ & $0.94$ & $1.0867e-2$ & $2.15$ & $4.7813e-1$ & $1.02$ & $4.7789e-1$ & $0.99$ \\ 
$1/40$    & $5.3561e-2$ & $0.97$ & $2.7257e-3$ & $1.99$ & $2.3968e-1$ & $1.00$ & $2.3953e-1$ & $1.00$ \\ 
$1/80$    & $2.7091e-2$ & $0.98$ & $7.2285e-4$ & $1.92$ & $1.1993e-1$ & $1.00$ & $1.1984e-1$ & $1.00$ \\ 
$1/160$   & $1.3624e-2$ & $0.99$ & $2.0194e-4$ & $1.84$ & $5.9976e-2$ & $1.00$ & $5.9932e-2$ & $1.00$ \\ 
\bottomrule
\end{tabular}
\end{center}
\end{table}


\begin{table}
\begin{center}
\caption{
Error and experimental convergence order for \eqref{steep_gradients} with $\mu = 10^{-2},\,\eta = 1$.}
\label{table:man_convergenceL2_m02}
\begin{tabular}{@{}ccccc|cccc@{}}
\toprule
 & \multicolumn{4}{c}{$\max_{1\leq n\leq \NT}\|U^n-u^n\|$} & \multicolumn{4}{c}{$\max_{1\leq n\leq \NT}\|C^n-c^n\|_1$} \\ \midrule
$h_0$ & LOW & Order & AFC & Order & LOW & Order & AFC & Order \\ \midrule
$1/10$    & $8.6443e-1$ &        & $1.2009e+0$ &        & $9.5335e-1$ &        & $9.5299e-1$ &       \\ 
$1/20$    & $6.0521e-1$ & $0.51$ & $1.5353e-1$ & $2.29$ & $4.8144e-1$ & $0.99$ & $4.7797e-1$ & $0.99$ \\ 
$1/40$    & $4.3452e-1$ & $0.48$ & $3.4673e-2$ & $2.15$ & $2.4184e-1$ & $1.00$ & $2.3954e-1$ & $1.00$ \\ 
$1/80$    & $2.8524e-1$ & $0.61$ & $8.8506e-3$ & $1.97$ & $1.2130e-1$ & $1.00$ & $1.1985e-1$ & $1.00$ \\ 
$1/160$   & $1.6840e-1$ & $0.76$ & $2.3454e-3$ & $1.91$ & $6.0788e-2$ & $1.00$ & $5.9933e-2$ & $1.00$ \\ 
\bottomrule
\end{tabular}
\end{center}
\end{table}

In Table \ref{table:man_convergenceL2_m00} we present the errors and the experimental orders of convergence in the $L^{2}$ and $H^{1}$ norms for \eqref{steep_gradients} with $\mu = 1$ at $T=1$, for both 
the low-order and AFC schemes. While both schemes achieve the optimal experimental order of convergence in the $H^{1}$ norm, this is not the case in the $L^{2}$ 
norm, where the AFC scheme attains nearly the optimal order, whereas the low-order scheme appears to converge linearly.

In addition, we set $\mu = 10^{-2}$, and in Table \ref{table:man_convergenceL2_m02}
we present the corresponding results for this case, where the AFC scheme seems again to be superior to low-order scheme.

\subsection{Convergence study of the stabilized schemes}
Next, we consider the following two sets of initial conditions for \eqref{Minimal_model_uc}.
The first one is
\begin{equation}\label{intial_cond_not_blow_up}
\begin{aligned}
u_0 & = 10\,e^{-10((x-0.5)^2 + (y-0.5)^2)} + 5,\quad 
c_0  = 0,
\end{aligned}
\end{equation}
where note that $\|u_0\|_{L^{1}} \approx 7.984 < 4\pi.$  
The  second set  is
\begin{equation}\label{intial_cond_not_blow_up2}
\begin{aligned}
u_0 & = \cos(2\pi x) + \cos(2\pi y) + 3,\quad 
c_0  =  \cos(2\pi x) + \sin(2\pi y) - 2\pi y + 9 ,
\end{aligned}
\end{equation}
where note that $\|u_0\|_{L^{1}} \approx 3 < 4\pi$. 
Thus, the solution $u$ in both examples does not blow-up. 
We consider again the non-stabilized scheme where we discretize in time \eqref{fem_uc} using the 
backward Euler method and the stabilized schemes given by 
\eqref{gen_fl_u_2D}. 
We assume a sequence of triangulations $\Th$ as described above with
$h_0 = 1/M$, $M = 10,20,40, 80$. The final time is chosen to be $T = 0.1,$ and we choose
$k = c_1\,h_0^{\log_{2} 3} \approx c_1\,h_0^{3/2}$ with $c_1 := \frac{2e-03}{10^{-\log_{2} 3}}.$
Since the exact solution of \eqref{Minimal_model_uc} is unknown, the underlying numerical reference solutions $U_{\mathrm{ref}}^{\mathrm{FEM}}\,, U_{\mathrm{ref}}^{\mathrm{LOW}}\,, U_{\mathrm{ref}}^{\mathrm{AFC}} \in \Sh$ at $t = T$ for each numerical scheme --i.e., FEM, LOW and AFC, respectively-- were obtained with $M = 320$ and a small time step $k = 5 \cdot 10^{-6},$ i.e., $\NT = 20000.$ Then, for different values of $\mu,$ the errors 
\begin{align*}
\|U_\sigma^T - U_{\mathrm{ref}}^\sigma\|,\,\,\|C_\sigma^T - C_{\mathrm{ref}}^\sigma\|_{1}, \quad \sigma \in \{\mathrm{FEM}, \mathrm{LOW}, \mathrm{AFC}\},
\end{align*}
are computed.

Initially, we consider different values of $\mu,$ such as $\mu = 1,\,10^{-2}$ for the system \eqref{Minimal_model_uc} by using the set of initial conditions given by \eqref{intial_cond_not_blow_up}. The results are presented in Tables \ref{table:exp_convergenceL2_m00}--\ref{table:exp_convergenceL2_m02}. While for $\mu=1,$ the results of all the numerical scheme are seems to be comparable, for $\mu = 10^{-2}$ the AFC scheme is superior to low-order scheme as it restores the optimal order in $L^{2}$ norm.


\begin{table}
\begin{center}
\caption{
Error and experimental convergence order for the initial conditions \eqref{intial_cond_not_blow_up} for \eqref{Minimal_model_uc} with $\mu = \eta = 1$.}
\label{table:exp_convergenceL2_m00}
\scalebox{0.85}{
\begin{tabular}{@{}ccccccc|cccccc@{}}
    \toprule
    & \multicolumn{6}{c}{$\|U_\sigma^T - U_{\mathrm{ref}}^\sigma\|$} & \multicolumn{6}{c}{$\|C_\sigma^T - C_{\mathrm{ref}}^\sigma\|_{1}$}\\ \midrule
    {$h_0$}  & {FEM}  	      & {Order}     &  {LOW}       	 & {Order}    &  {AFC}        & {Order}    & {FEM}           & {Order}     &  {LOW}         & {Order}    &  {AFC}        & {Order} \\ \hline
    {$1/10$} & {$2.175e-2$}    &            & {$2.063e-2$}   &           & {$2.489e-2$}  &           & {$9.515e-3$}    &            & {$1.057e-2$}   &           & {$1.151e-2$}  &   \\ \hline
    {$1/20$} & {$5.340e-3$}    & {$2.03$}   & {$4.978e-3$}   & {$2.05$}  & {$6.692e-3$}  & {$1.89$}  & {$4.259e-3$}    & {$1.16$}   & {$4.380e-3$}   & {$1.27$}  & {$4.617e-3$}  & {$1.32$}\\ \hline
    {$1/40$} & {$1.266e-3$}    & {$2.08$}   & {$1.134e-3$}   & {$2.13$}  & {$1.815e-3$}  & {$1.88$}  & {$2.049e-3$}    & {$1.06$}   & {$2.060e-3$}   & {$1.09$}  & {$2.111e-3$}  & {$1.13$}\\ \hline
    {$1/80$} & {$2.936e-4$}    & {$2.11$}   & {$2.351e-4$}   & {$2.27$}  & {$4.851e-4$}  & {$1.90$}  & {$9.925e-4$}    & {$1.05$}   & {$9.924e-4$}   & {$1.05$}  & {$1.002e-3$} & {$1.07$}\\ \bottomrule
\end{tabular}
}
\end{center}
\end{table}


\begin{table}
\begin{center}
\caption{
Error and experimental convergence order for the initial conditions \eqref{intial_cond_not_blow_up} for \eqref{Minimal_model_uc} with $\mu = 10^{-2},\,\eta = 1$.}
\label{table:exp_convergenceL2_m02}
\scalebox{0.85}{
\begin{tabular}{@{}ccccccc|cccccc@{}}
    \toprule
    & \multicolumn{6}{c}{$\|U_\sigma^T - U_{\mathrm{ref}}^\sigma\|$} & \multicolumn{6}{c}{$\|C_\sigma^T - C_{\mathrm{ref}}^\sigma\|_{1}$}\\ \midrule
    {$h_0$}  & { FEM}          & {Order}     &  {LOW}         & {Order}     &  {AFC}         & {Order}     & {FEM}           & {Order}     &  {LOW}         & {Order}     &  {AFC}        & {Order} \\ \hline
    {$1/10$} & {$3.968e-1$}    &            & {$1.259e+0$}   &            & {$1.018e+0$}   &            & {$1.374e-1$}    &            & {$1.400e-1$}   &            & {$1.392e-1$}  &     \\ \hline
    {$1/20$} & {$1.077e-1$}    & {$1.88$}   & {$6.139e-1$}   & {$1.04$}   & {$3.650e-1$}   & {$1.48$}   & {$6.998e-2$}    & {$0.97$}   & {$7.342e-2$}   & {$0.93$}   & {$7.034e-2$}  & {$0.98$}\\ \hline
    {$1/40$} & {$2.541e-2$}    & {$2.08$}   & {$2.731e-1$}   & {$1.17$}   & {$8.650e-2$}   & {$2.08$}   & {$3.494e-2$}    & {$1.00$}   & {$3.714e-2$}   & {$0.98$}   & {$3.498e-2$}  & {$1.01$}\\ \hline
    {$1/80$} & {$6.742e-3$}    & {$1.91$}   & {$1.128e-1$}   & {$1.28$}   & {$1.742e-2$}   & {$2.31$}   & {$1.706e-2$}    & {$1.03$}   & {$1.804e-2$}   & {$1.04$}   & {$1.707e-2$}  & {$1.04$}\\    \bottomrule
\end{tabular}
}
\end{center}
\end{table}

Next, we test all the numerical schemes for the set of initial conditions given by \eqref{intial_cond_not_blow_up2}. Unlike the previous set of initial conditions, the superiority of the AFC scheme over the low-order scheme can also be observed when $\mu = 1$. In particular, Table \ref{table:tri_convergenceL2_m00} clearly shows that the low-order scheme exhibits only first-order convergence $\mathcal{O}(h)$ in the $L^{2}$ norm, while both the AFC scheme and the standard FEM achieve the optimal order of convergence in the $L^{2}$ norm. Furthermore, in the $H^{1}$ norm, all the schemes exhibit the same order of convergence, but the absolute errors of the AFC scheme and the standard FEM are smaller than those of the low-order scheme.

It easily verifiable that the experimental order of convergence in $L^{2}$ for the AFC method is optimal,
 in contrast to the theoretical order of convergence proved on Theorem \ref{theorem:error_estimates_low_order_2D_main}. 
 The sub-optimal result of Theorem \ref{theorem:error_estimates_low_order_2D_main} comes from the estimation of the chemotactic term, see, e.g., Lemma \ref{lemma:chemotaxis_bounds_2D-1}.


\begin{table}
\begin{center}
\caption{Error and experimental convergence order for the initial conditions \eqref{intial_cond_not_blow_up2} for \eqref{Minimal_model_uc} with $\mu = \eta = 1$.}
\label{table:tri_convergenceL2_m00}
\scalebox{0.85}{
\begin{tabular}{@{}ccccccc|cccccc@{}}
    \toprule
    & \multicolumn{6}{c}{$\|U_\sigma^T - U_{\mathrm{ref}}^\sigma\|$} & \multicolumn{6}{c}{$\|C_\sigma^T - C_{\mathrm{ref}}^\sigma\|_{1}$}\\ \midrule
    {$h_0$}  & {FEM}           & {Order}     &  {LOW}         & {Order}    &  {AFC}         & {Order}    & { FEM}          & {Order}     &  {LOW}         & {Order}    &  {AFC}         & {Order} \\ \hline
    {$1/10$} & {$4.681e-2$}    &            & {$2.665e-1$}   &           & {$5.155e-2$}   &           & {$2.617e-1$}    &            & {$2.630e-1$}   &           & {$2.687e-1$}   &    \\ \hline
    {$1/20$} & {$1.224e-2$}    & {$1.94$}   & {$1.341e-1$}   & {$0.99$}  & {$1.436e-2$}   & {$1.84$}  & {$1.312e-1$}    & {$1.00$}   & {$1.317e-1$}   & {$1.00$}  & {$1.324e-1$}   & {$1.02$}\\	\hline
    {$1/40$} & {$3.095e-3$}    & {$1.98$}   & {$6.401e-2$}   & {$1.07$}  & {$4.183e-3$}   & {$1.78$}  & {$6.534e-2$}    & {$1.01$}   & {$6.568e-2$}   & {$1.00$}  & {$6.555e-2$}   & {$1.01$}\\ \hline
    {$1/80$} & {$7.591e-4$}    & {$2.03$}   & {$2.777e-2$}   & {$1.20$}  & {$1.212e-3$}   & {$1.79$}  & {$3.190e-2$}    & {$1.03$}   & {$3.206e-2$}   & {$1.03$}  & {$3.193e-2$}   & {$1.04$}\\    \bottomrule
\end{tabular}
}
\end{center}
\end{table}

\section{Conclusions}
In this paper, we presented a finite element error analysis for a stabilized finite element scheme for the Chemotaxis systems based on the AFC method. In particular, we discretized space using continuous piecewise linear finite elements and time by the backward Euler method. Under assumptions for the triangulation used for the space discretization and the size of the time step $k$, we showed  that the resulting coupled non-linear scheme have a unique solution which remain non-negative. Further, we showed that the resulting discrete solution remains bounded and derived error estimates in $L^{2}$ and $H^{1}-$norm in space. Numerical experiments in two dimensions were presented for both the standard FEM and stabilized schemes. In the numerical experiments, we computed the experimental order of convergence and we observed optimal order of convergence in some cases for the AFC method in $L^{2}-$norm, while for low-order, the order is linear.

\section*{Acknowledgments}
\textit{Acknowledgments.} The research of C. Pervolianakis was supported by the Hellenic Foundation for Research and Innovation (HFRI) under the HFRI PhD Fellowship grant (Fellowship Number: 837).

\bibliographystyle{plain} 
\bibliography{ref} 

\end{document}